# Efficient Multi-Accuracy Computations of Complex Functions with Complex Arguments

MOFREH R. ZAGHLOUL, United Arab Emirates University[1]

We present an efficient multi-accuracy algorithm for the computations of a set of special functions of a complex argument, $z=x+iy$. These functions include the complex probability function $w(z)$, and closely related functions such as the error function $erf(z)$, complementary error function $erfc(z)$, imaginary error function $erfi(z)$, scaled complementary error function, $erfcx(z)$, the plasma dispersion function $Z(z)$, Dawson's function $Daw(z)$, and Fresnel integrals $S(z)$ and $C(z)$. Computational results from the present algorithm are compared with results from competitive algorithms and widely used software packages showing superior accuracy and efficiency of the present algorithm. In particular, the present results highlight concerns about the accuracy of evaluating such special functions using commercial packages like *Mathematica* and free/open source packages like the *MIT C++* package.

---

Author's addresses: M. Zaghloul, Department of Physics, College of Sciences, United Arab Emirates University, Al-Ain, 15551, UAE.

## 1. INTRODUCTION

The evaluation of a group of functions or integrals of complex arguments (of central importance to many fields of physics and applied sciences) may depend in one way or another on the accurate evaluation of the complex probability function (CPF), also known as the Faddeyeva (or Faddeeva) function, $w(z)$. Accordingly, accurate and efficient evaluation of the complex probability function is a topic that continues to gain extensive interest in the literature. The function can be written in different forms as,

$$w(z) = e^{(-iz)^2} erfc(-iz) = e^{-z^2}(1 + i\, erfi(z)) = erfcx(-iz) \tag{1},$$

where $erfc(z)$ is the complementary error function, $erfi(z)$ is the imaginary error function, and the designation "$erfcx$" is also used in the literature to refer to the scaled complementary error function for complex and real arguments. The real and imaginary parts of the function are known as the real and imaginary Voigt functions.

Many studies and algorithms have been introduced in the literature for the analysis and evaluation of the function [Faddeyeva and Terent'ev 1961, Young C. 1965; Armstrong 1967; Gautschi 1969; 1970; Hui et al. 1978; Humlíček 1982; Dominguez et al. 1987; Poppe and Wijers. 1990a,b; Lether and Wenston 1991; Schreier 1992; Shippony and Read 1993; Weideman 1994; Wells 1999; Luque et al. 2005; Letchworth and Benner 2007; Abrarov et al. 2010a,b; Zaghloul and Ali 2011, Johnson 2012; Boyer and Lynas-Gray 2014, Zaghloul 2015, Zaghloul 2016, Zaghloul 2017]. Algorithm 916 [Zaghloul and Ali 2011] is one of these algorithms, originally introduced as a *Matlab* function "*Faddeyeva.m*" for multi-accuracy evaluation of the function, and is recognized for its remarkable accuracy. Efficiency improvements and a *Fortran* translation of the algorithm were introduced [Zaghloul 2016] to increase practical reliability of the algorithm. The latter version of Algorithm 916 (*Faddeyeva_v2*) is both accurate and efficient.

In the present work, we present a more efficient algorithm for multi-accuracy computations of the function and then use it as the basis of the computation of a set of related special functions of complex arguments.

The foundation and details of the present algorithm for calculating the complex probability function, are described in Section 2 followed by efficiency benchmark and speed comparison in Section 3. A summary of mathematical relations and the methods used for computing the set of related special functions of complex arguments is given in Section 4, while accuracy verification and comparison are presented in Section 5.

## 2. Algorithm Foundation

In this section we introduce the foundation of the method used herein for a more efficient multi-accuracy evaluation of the complex probability function, CPF, compared to other competitive algorithms in the literature including the most recent version of Algorithm 916, *Faddeyeva_v2*, and the free/open-source "*Faddeeva Package*" developed in *C++* at the *Massachusetts Institute of Technology* (*MIT*) [Johnson 2012].

**2.1. Accurate Asymptotic Representation**

For large values of $|z|$, the present algorithm makes use of two asymptotic expressions for the Faddeyeva function:

(a) the asymptotic approximation using Laplace continued fraction [Faddeyeva and Terent'ev 1961, Abramowitz and Stegun 1964, Gautschi 1970],

$$w(z) = \frac{i}{\sqrt{\pi}} \frac{1}{z-} \frac{\frac{1}{2}}{z-} \frac{1}{z-} \frac{\frac{3}{2}}{z-} \frac{2}{z-} \ldots\ldots\ldots, \qquad y > 0 \qquad (2)$$

where the continued fraction needs to be truncated at some convergent for practical evaluation. For large values of $|z|$, a few convergents may be used to obtain a desired accuracy. Table 1 summarizes the rational approximations resulting from the first six convergents of the Laplace continued fraction along with their regions of applicability for different targeted accuracies, as drawn from a systematic accuracy check using *Algorithm 916* as a reference;

(b) another useful series approximation for the Faddeyeva function as $z \to \infty$ and $|\arg z| < 3\pi/4$ [Abramowitz and Stegun, 1964 7.1.23 & Zwillinger 2003, 6.13.3] can be written as

$$w(z) \sim \frac{i}{z\sqrt{\pi}} \sum_{n=0}^{\infty} \frac{(2n)!}{n!\, 2^n\, (2z^2)^n} \approx \frac{i}{z\sqrt{\pi}} \sum_{m=0}^{\infty} \frac{(2m-1)!!}{(2z^2)^m}$$
$$\sim \frac{i}{z\sqrt{\pi}} \left( 1 + \frac{1}{2z^2} + \frac{3}{4z^4} + \frac{15}{8z^6} + \frac{105}{16z^8} + \frac{945}{32z^{10}} + \frac{10395}{64z^{12}} + \frac{135135}{128z^{14}} + \frac{2027025}{256z^{16}} + .. \right) \qquad (3),$$

where $(2m+1)!! = 1 \cdot 3 \cdot 5 \cdots (2m+1)$ and $(-1)!! = 1$. The above series can be calculated efficiently in the form

$$w(z) \sim \frac{i}{z\sqrt{\pi}} \left( 1 + \alpha \left( 1 + \alpha \left( 3 + \alpha \left( 15 + \alpha \left( 105 + \alpha \left( 945 + \alpha \left( 10395 + \alpha \left( 135135 + \ldots \right)\right)\right)\right)\right)\right)\right)\right) \quad (4),$$

with $\alpha = 1/(2z^2)$.

A few to several terms of this series may be used to obtain the desired accuracy for $|z| \gg 1$. Yet again, *Algorithm 916* has been used as a reference for a process of systematic accuracy check and Table 2 summarizes the terms retained from the series and the corresponding regions of applicability (in the first quarter) for different targeted accuracies.

Some interesting features and useful findings can be derived from the results in Tables 1 and 2:

(i) Using six convergents from Laplace continued fraction or retaining six terms from the series in Eq. (3) is sufficient to calculate the Faddeyeva function with an accuracy up 13 digits after the decimal point for the regions $|z|^2 \geq 400$ and $|z|^2 \geq 277$, respectively. These ranges can be extended to smaller values of $|z|$ for lower accuracies as appears from the Tables. This necessarily covers the region $|z| \geq \sqrt{-ln(R_{min})}$, where $R_{min}$ is the smallest positive normalized floating point number, which means that one can get rid of one of the loops used in the recent version of Algorithm 916 [Zaghloul 2016] for the calculation of the function when $|x| \geq \sqrt{-ln(R_{min})}$ and replace it either by convergents from the Laplace continued fraction or approximations from the series in Eq. (3). Such modification saves considerable execution time.

(ii) Using more than six terms of the series (3) it is possible to extend the domain of application of the approximation to $|z|^2 \geq 127$ for 13 digits accuracy and down to $|z|^2 \geq 120$ for 10 digits accuracy, which covers a considerable part of the domain of the remaining computational loop used for the calculation of the function in the recent version of Algorithm 916 [Zaghloul 2016] for $|x| < \sqrt{-ln(R_{min})}$ with additional saving in the execution time compared to "*Faddeyeva_v2*".

(iii) Increasing the number of convergents from the Laplace continued fraction to more than 6 or the number of terms retained from the series (3) to more than 9 terms does not extend the region of applicability or improve the efficiency and as a result they are not shown in Tables 1 and 2.

(iv) As the computational burden increases with the number of convergents and/or the number of terms retained from the series (3), a modest effort is devoted to choose the efficient expressions used for obtaining the targeted accuracy as explained in the algorithm components given below.

(v) As a derivative, Table 1 reveals the reason for the loss of accuracy experienced with Algorithm 680 [Gautschi 1970], where Laplace continued fraction was used down to $x = 6.3$ which is outside the range of applicability for small values of $y$, for any of the targeted accuracies from 4 to 13 significant digits as shown in Table 1. , see [Zaghloul 2019] for more details).

**Table 1:** Rational approximations from Laplace continued fraction vs regions of applicability for different targeted accuracies

| # of Convergents | Rational approximation of $w(z)$ | Accuracy | | | | | | | | | |
|---|---|---|---|---|---|---|---|---|---|---|---|
| | | $\varepsilon=10^{-13}$ | $\varepsilon=10^{-12}$ | $\varepsilon=10^{-11}$ | $\varepsilon=10^{-10}$ | $\varepsilon=10^{-9}$ | $\varepsilon=10^{-8}$ | $\varepsilon=10^{-7}$ | $\varepsilon=10^{-6}$ | $\varepsilon=10^{-5}$ | $\varepsilon=10^{-4}$ |
| | | Region of Applicability $\|z\|^2=(x^2+y^2) \geq$ | | | | | | | | | |
| 1 | $\approx \dfrac{i}{z\sqrt{\pi}}$ | $1.5\times10^{13}$ | $1.2\times10^{12}$ | $1.5\times10^{11}$ | $1.5\times10^{10}$ | $1.4\times10^{9}$ | $1.3\times10^{8}$ | $1.5\times10^{7}$ | $1.451\times10^{6}$ | $1.5\times10^{5}$ | **16000** |
| 2 | $\approx \dfrac{iz}{\sqrt{\pi}\,(z^2-0.5)}$ | $10^{8}$ | $1.9\times10^{6}$ | $5.0\times10^{5}$ | $2.0\times10^{5}$ | 50000 | 16000 | 5000 | 1600 | 510 | **160** |
| 3 | $\approx \dfrac{i}{z\sqrt{\pi}}\dfrac{(z^2-1)}{(z^2-1.5)}$ | $3.8\times10^{4}$ | 17500 | 8100 | 3750 | 1750 | 810 | 380 | 180 | 110 | **107** |
| 4 | $\approx \dfrac{iz}{\sqrt{\pi}}\dfrac{(z^2-2.5)}{z^2(z^2-3)+0.75}$ | $3.5\times10^{3}$ | 1950 | 1085 | 610 | 345 | 195 | 115 | **111** | **109** | 107 |
| 5 | $\approx \dfrac{i}{z\sqrt{\pi}}\dfrac{z^2(z^2-4.5)+2.0}{z^2(z^2-5)+3.75}$ | 12000 | 550 | 340 | 215 | 137 | **116** | **114** | 111 | 109 | 107 |
| 6 | $\approx \dfrac{iz}{\sqrt{\pi}}\dfrac{(z^2(z^2-7)+8.25)}{z^2(z^2(z^2-7.5)+11.25)-1.875}$ | **400** | **235** | **162** | **122** | **118** | 116 | 114 | 111 | 109 | 107 |

**Table 2**: Series approximation of the Faddeyeva function Eq. (3): Terms retained vs region of applicability for different targeted accuracies

| m | Approximation of $w(z) \approx \frac{i}{z\sqrt{\pi}}(1+f)$ | Accuracy | | | | | | | | | |
|---|---|---|---|---|---|---|---|---|---|---|---|
| | | $10^{-13}$ | $10^{-12}$ | $10^{-11}$ | $10^{-10}$ | $10^{-9}$ | $10^{-8}$ | $10^{-7}$ | $10^{-6}$ | $10^{-5}$ | $10^{-4}$ |
| | | Region of Applicability $\|z\|^2=(x^2+y^2) \geq$ | | | | | | | | | |
| 0 | $f \approx 0$ | $1.5\times10^{13}$ | $1.2\times10^{12}$ | $1.5\times10^{11}$ | $1.5\times10^{10}$ | $1.4\times10^{9}$ | $1.3\times10^{8}$ | $1.5\times10^{7}$ | $1.451\times10^{6}$ | $1.5\times10^{5}$ | **16000** |
| 1 | $f \approx \alpha\,(1)$ | $5.0\times10^{6}$ | $1.74\times10^{6}$ | $5.05\times10^{5}$ | $1.71\times10^{5}$ | 45626 | 12136 | 6096 | 1937 | 606 | **191** |
| 2 | $f \approx (\alpha\,(1+\alpha\,(3)))$ | 21000 | 14360 | 10245 | 5079 | 2358 | 1096 | 509 | 236 | 112 | **107** |
| 3 | $f \approx (\alpha\,(1+\alpha\,(3+\alpha\,(15))))$ | 6800 | 2720 | 1554 | 877 | 494 | 278 | 157 | 112 | **109** | 107 |
| 4 | $f \approx (\alpha(1+\alpha(3+\alpha(15+\alpha(105)))))$ | 2500 | 799 | 505 | 319 | 201 | 127 | **114** | **111** | 109 | 107 |
| 5 | $f \approx (\alpha(1+\alpha(3+\alpha(15+\alpha(105+\alpha(945))))))$ | 750 | 360 | 356 | 167 | 119 | **116** | 114 | 111 | 109 | 107 |
| 6 | $f \approx (\alpha(1+\alpha(3+\alpha(15+\alpha(105+\alpha(945+\alpha(110395)))))))$ | 277 | 208 | 149 | 122 | **118** | 116 | 114 | 111 | 109 | 107 |
| 7 | $f \approx (\alpha(1+\alpha(3+\alpha(\ldots\ldots+\alpha(110395+\alpha(135135)))))))$ | 177 | 140 | **123** | **120** | 118 | 116 | 114 | 111 | 109 | 107 |
| 8 | $f \approx (\alpha(1+\alpha(3+\alpha(\ldots\ldots+\alpha(135135+\alpha(2027025))))))))$ | 129 | **125** | 123 | 120 | 118 | 116 | 114 | 111 | 109 | 107 |
| 9 | $f \approx (\alpha(1+\alpha(3+\alpha(\ldots\ldots+\alpha(2027025+\alpha(34459425)))))))))$ | **127** | 125 | 123 | 120 | 118 | 116 | 114 | 111 | 109 | 107 |

## 2.2. Dawson's Function of a Real Argument

Systematic numerical experiments using Algorithm 916 as a reference showed that the region of applicability of the convergents from the Laplace continued fraction given in Table 1 and the terms from the series (3) given in Table 2 can be extended further to smaller values of $|z|$ except for very small values of $y$. For such regions of very small values of $y$, Algorithm 916 can be used, however, it was found that using a few terms from the Taylor expansion of the Dawson's integral, $Daw(z)$, to approximate $w(z)$ where $w(z)=exp(-z^2)+2i\,Daw(z)/\pi^{1/2}$, is more efficient.

Expanding the Dawson's integral function about some point $x$ on the $x$-axis produces a recursion relation among the series coefficients [Armstrong 1967, Shippony and Read 1993],

$$d_0 = Daw(x), \quad d_1 = 1 - 2xd_0, \quad d_{n+1} = -\frac{2}{n+1}(xd_n + d_{n-1}) \text{ for } n = 1,2,\ldots \quad (5)$$

Dawson's function of a real argument, $Daw(x)$, can be calculated using Algorithm 715 [Cody 1993] which can also be used for the calculation of $Daw(z)$ when $y=0$.

## 2.3. Algorithm Components

In the present algorithm, the first quadrant in the computational domain is divided into six main regions ($I$ to $VI$) as shown in Figure 1. For efficiency reasons, regions $IV$ and $V$ may be further split into two or more sub-regions as indicated, for region $IV$, by the dotted borders in the figure. For efficient computations, the borders for the main regions and sub-regions vary according to the targeted accuracy.

The algorithm considers 10 different targeted accuracies between 4 to 13 digits after the decimal point. For all targeted accuracies, the algorithm uses the first three convergents from Laplace continued fraction to approximate the function asymptotically for the outer three regions of large $|z|$ (1 convergent for region $I$, 2 convergents for region $II$, and 3 convergents for region $III$) though with different borders as shown in the Tables 3a-j below. Because the asymptotic expressions from Laplace continued fractions are not valid on the real axis ($y=0$), the Faddeyeva function is calculated on the real axis using Dawson's integral, $Daw(x)$, where

$$w(x) = \exp(-x^2) + \frac{2i}{\sqrt{\pi}} Daw(x) \quad (6)$$

This expression has been also implemented in the loop of Algorithm 916 (whenever used herein) to replace the calculation of the function on the real axis using the approximation $|y| \to R_{min}$ when $y=0$.

To a great extent the details of the algorithm for 4 significant-digits accuracy is similar to Algorithm 985 [Zaghloul 2017] with slight change in the borders of regions and the calculation of $w(x)$ using Eq. (6) for the outer three regions. However, for five significant figures, in addition to new different borders, there is a region where none of the methods used in the latter can secure the targeted accuracy of five significant digits.

For this particular region, the loop used for the region $|x| < \sqrt{-ln(R_{min})}$ in Algorithm 916 is used to calculate the function to the targeted accuracy as explained in Table 3-b.

The use of a few terms from the Taylor expansion of the Dawson's integral to approximate the function for very small values of *y* (sub-region of *V*) is implemented for accuracies from 6 to 13 significant digits as shown in Tables 3c-j below for its favorable efficiency. Although omitted from Tables 3c-j, for brevity, each region in which Taylor expansion of the Dawson's integral is used to approximate the function is divided into finer sub-regions using different number of terms for efficient implementation. In addition, when the value of the imaginary parameter *y* is very small such that one or two terms of the Taylor series are sufficient to secure the targeted accuracy, one can further save execution time by approximating the term exp(-$z^2$) by exp($y^2$-$x^2$) (1.0-*2i xy*) avoiding calling the *sine* and *cosine* intrinsic functions. When the approximation for region *IV* in Humlíček's w4 algorithm is used to approximate the function for very small values of *y*, a similar implementation or even a simpler one (exp(-$z^2$) ≈ exp(-$x^2$)) is used for the same reason. This is clearly understandable in the light of the low accuracy sought in such cases (only four or five significant-digits cases) for which Humlíček's expression is called for a very narrow region of the computational domain.

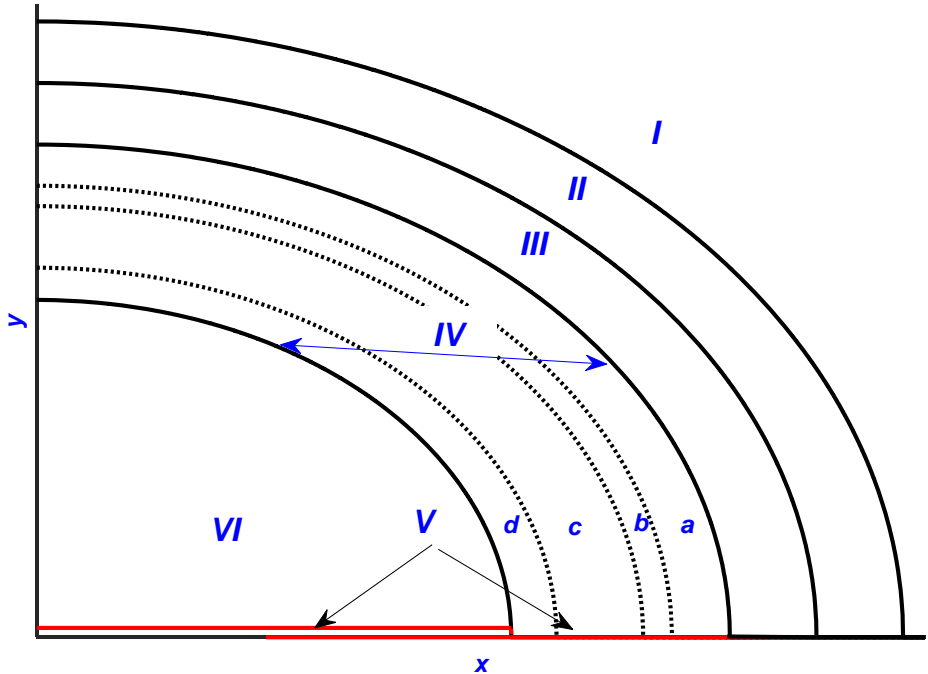

**Figure 1.** A schematic diagram for the regions of the computational domain used in the present algorithm

**Table 3-a. Details of the present algorithm component for four significant digits**

| 4 Significant figures | | |
|---|---|---|
| Region | Borders | Method |
| I | $\|z\|^2 \geq 16000.0$ | 1 convergent (Laplace continued fractions) |
| II | $16000.0 > \|z\|^2 \geq 160.0$ | 2 convergents (Laplace continued fractions) |
| III | $160.0 > \|z\|^2 \geq 107.0$ | 3 convergents (Laplace continued fractions) |
| IV | $107.0 > \|z\|^2 \geq 28.5$ & $y^2 \geq 6 \times 10^{-14}$ | 4 convergents (Laplace continued fractions) |
| V | $107.0 > \|z\|^2 \geq 28.5$ & $y^2 < 6 \times 10^{-14}$ $\cup$ $28.5 > \|z\|^2 \geq 3.5$ & $y^2 < 0.026$ | Approx. for region IV in Humlíček's w4 algorithm |
| VI | Otherwise | Hui's p6 approximation |

**Table 3-b. Details of the present algorithm component for five significant digits**

| 5 Significant figures | | |
|---|---|---|
| Region | Borders | Method |
| I | $\|z\|^2 \geq 150000.0$ | 1 convergent (Laplace continued fractions) |
| II | $150000.0 > \|z\|^2 \geq 510.0$ | 2 convergents (Laplace continued fractions) |
| III | $510.0 > \|z\|^2 \geq 110.0$ | 3 convergents (Laplace continued fractions) |
| IV | $110.0 > \|z\|^2 \geq 109.0$ $\cup$ $109.0 > \|z\|^2 \geq 39.0$ & $y^2 \geq 10^{-9}$ | 4 convergents (Laplace continued fractions) |
| V | $109.0 > \|z\|^2$ & $y^2 < 10^{-9}$ | Approx. for region IV in Humlíček's w4 algorithm |
|  | $109.0 > \|z\|^2$ & $0.1 \geq y^2 > 10^{-9}$ | $Daw715(x)$ and a few terms of Taylor series |
|  | $39.0 > \|z\|^2$ & $0.27 \geq y^2 > 0.1$ | Algorithm 916 (loop used for $x < (-\ln R_{min})^{1/2}$) |
| VI | $39.0 > \|z\|^2$ & $y^2 > 0.27$ | Hui's p6 Approximation |

**Table 3-c. Details of the present algorithm component for six significant digits**

| 6 Significant figures | | |
|---|---|---|
| Region | Borders | Method |
| I | $\|z\|^2 \geq 1451000.0$ | 1 convergent (Laplace continued fractions) |
| II | $1451000.0 > \|z\|^2 \geq 1600.0$ | 2 convergents (Laplace continued fractions) |
| III | $1600.0 > \|z\|^2 \geq 180.0$ | 3 convergents (Laplace continued fractions) |
| IV | $180.0 > \|z\|^2 \geq 111.0$ | 4 convergents (Laplace continued fractions) |
| V | $111.0 > \|z\|^2$ & $y^2 \leq 10^{-2}$ | $Daw715(x)$ and a few terms of Taylor series |
|  | $111.0 > \|z\|^2$ & $1 > y^2 > 10^{-2}$ | Algorithm 916 (loop used for $x < (-\ln R_{min})^{1/2}$) |
| VI | $111.0 > \|z\|^2$ & $y^2 \geq 1$ | Hui's p6 Approximation |

**Table 3-d. Details of the present algorithm component for seven significant digits**

| 7 Significant figures | | |
|---|---|---|
| Region | Borders | Method |
| I | $\|z\|^2 \geq 1.5 \times 10^7$ | 1 convergent (Laplace continued fractions) |
| II | $1.5 \times 10^7 > \|z\|^2 \geq 5010.0$ | 2 convergents (Laplace continued fractions) |
| III | $5010.0 > \|z\|^2 \geq 380.0$ | 3 convergents (Laplace continued fractions) |
| IV | $380.0 > \|z\|^2 \geq 115.0$ | 4 convergents (Laplace continued fractions) |
|  | $115.0 > \|z\|^2 \geq 114.0$ | 5 convergents (Laplace continued fractions) |
| V | $114.0 > \|z\|^2$ & $y^2 \leq 10^{-2}$ | $Daw715(x)$ and a few terms of Taylor series |
| VI | Otherwise | Algorithm 916 (loop used for $x < (-\ln R_{min})^{1/2}$) |

**Table 3-e. Details of the present algorithm component for eight significant digits**

| Region | Borders | Method |
|---|---|---|
| \multicolumn{3}{c}{8 Significant figures} | | |
| I | $\|z\|^2 \geq 1.3 \times 10^8$ | *1 convergent (Laplace continued fractions)* |
| II | $1.3 \times 10^8 > \|z\|^2 \geq 16000.0$ | *2 convergents (Laplace continued fractions)* |
| III | $16000.0 > \|z\|^2 \geq 810.0$ | *3 convergents (Laplace continued fractions)* |
| IV | $810.0 > \|z\|^2 \geq 195.0$ | *4 convergents (Laplace continued fractions)* |
|  | $195.0 > \|z\|^2 \geq 116.0$ | *5 convergents (Laplace continued fractions)* |
| V | $116.0 > \|z\|^2$ & $y^2 \leq 10^{-3}$ | *Daw715(x) and a few terms of Taylor series* |
| VI | *Otherwise* | *Algorithm 916 (loop used for $x < (-\ln R_{min})^{1/2}$)* |

**Table 3.f. Details of the present algorithm component for nine significant digits**

| Region | Borders | Method |
|---|---|---|
| \multicolumn{3}{c}{9 Significant figures} | | |
| I | $\|z\|^2 \geq 1.4 \times 10^9$ | *1 convergent (Laplace continued fractions)* |
| II | $1.4 \times 10^9 > \|z\|^2 \geq 50000.0$ | *2 convergents (Laplace continued fractions)* |
| III | $50000.0 > \|z\|^2 \geq 1750.0$ | *3 convergents (Laplace continued fractions)* |
| IV | $1750.0 > \|z\|^2 \geq 345.0$ | *4 convergents (Laplace continued fractions)* |
|  | $345.0 > \|z\|^2 \geq 137.0$ | *5 convergents (Laplace continued fractions)* |
|  | $137.0 > \|z\|^2 \geq 118.0$ | *6 convergents (Laplace continued fractions)* |
| VI | $118.0 > \|z\|^2$ & $y^2 \leq 10^{-3}$ | *Daw715(x) and a few terms of Taylor series* |
| VI | *Otherwise* | *Algorithm 916 (loop used for $x < (-\ln R_{min})^{1/2}$)* |

**Table 3.g. Details of the present algorithm component for 10 significant digits**

| Region | Borders | Method |
|---|---|---|
| \multicolumn{3}{c}{10 Significant figures} | | |
| I | $\|z\|^2 \geq 1.5 \times 10^{10}$ | *1 convergent (Laplace continued fractions)* |
| II | $1.5 \times 10^{10} > \|z\|^2 \geq 200000.0$ | *2 convergents (Laplace continued fractions)* |
| III | $200000.0 > \|z\|^2 \geq 3750.0$ | *3 convergents (Laplace continued fractions)* |
| IV | $3750.0 > \|z\|^2 \geq 610.0$ | *4 convergents (Laplace continued fractions)* |
|  | $610.0 > \|z\|^2 \geq 215.0$ | *5 convergents (Laplace continued fractions)* |
|  | $215.0 > \|z\|^2 \geq 122.0$ | *6 convergents (Laplace continued fractions)* |
|  | $122.0 > \|z\|^2 \geq 120.0$ | *6 terms of the series in (4)* |
| V | $120.0 > \|z\|^2$ & $y^2 \leq 10^{-3}$ | *Daw715(x) and a few terms of Taylor series* |
| VI | *Otherwise* | *Algorithm 916 (loop used for $x < (-\ln R_{min})^{1/2}$)* |

**Table 3.h. Details of the present algorithm component for 11 significant digits**

| Region | Borders | Method |
|---|---|---|
| \multicolumn{3}{c}{11 Significant figures} | | |
| I | $\|z\|^2 \geq 1.5 \times 10^{11}$ | *1 convergent (Laplace continued fractions)* |
| II | $1.5 \times 10^{11} > \|z\|^2 \geq 500000.0$ | *2 convergents (Laplace continued fractions)* |
| III | $500000.0 > \|z\|^2 \geq 8100.0$ | *3 convergents (Laplace continued fractions)* |
| IV | $8100.0 > \|z\|^2 \geq 1085.0$ | *4 convergents (Laplace continued fractions)* |
|  | $1085.0 > \|z\|^2 \geq 340.0$ | *5 convergents (Laplace continued fractions)* |
|  | $340.0 > \|z\|^2 \geq 162.0$ | *6 convergents (Laplace continued fractions)* |
|  | $162.0 > \|z\|^2 \geq 123.0$ | *7 terms of the series in (4)* |
| V | $123.0 > \|z\|^2$ & $y^2 \leq 10^{-3}$ | *Daw715(x) and a few terms of Taylor series* |
| VI | *Otherwise* | *Algorithm 916 (loop used for $x < (-\ln R_{min})^{1/2}$)* |

**Table 3.i. Details of the present algorithm component for 12 significant digits**

| Region | Borders | Method |
|---|---|---|
| 12 Significant figures | | |
| I | $|z|^2 \geq 1.2 \times 10^{12}$ | *1 convergent (Laplace continued fractions)* |
| II | $1.2 \times 10^{12} > |z|^2 \geq 1900000.0$ | *2 convergents (Laplace continued fractions)* |
| III | $1900000.0 > |z|^2 \geq 17500.0$ | *3 convergents (Laplace continued fractions)* |
| IV | $17500.0 > |z|^2 \geq 1950.0$ | *4 convergents (Laplace continued fractions)* |
|  | $1950.0 > |z|^2 \geq 550.0$ | *5 convergents (Laplace continued fractions)* |
|  | $550.0 > |z|^2 \geq 235.0$ | *6 convergents (Laplace continued fractions)* |
|  | $235.0 > |z|^2 \geq 125.0$ | *8 terms of the series in (4)* |
| V | $125.0 > |z|^2$ & $y^2 \leq 10^{-3}$ | *Daw715(x) and a few terms of Taylor series* |
| VI | *Otherwise* | *Algorithm 916 (loop used for $x < (-\ln R_{min})^{1/2}$)* |

**Table 3.j. Details of the present algorithm component for 13 significant digits**

| Region | Borders | Method |
|---|---|---|
| 13 Significant figures | | |
| I | $|z|^2 \geq 1.5 \times 10^{13}$ | *1 convergent (Laplace continued fractions)* |
| II | $1.5 \times 10^{13} > |z|^2 \geq 1.0 \times 10^8$ | *2 convergents (Laplace continued fractions)* |
| III | $1.0 \times 10^8 > |z|^2 \geq 38000.0$ | *3 convergents (Laplace continued fractions)* |
| IV | $38000.0 > |z|^2 \geq 3500.0$ | *4 convergents (Laplace continued fractions)* |
|  | $3500.0 > |z|^2 \geq 1200.0$ | *5 convergents (Laplace continued fractions)* |
|  | $1200.0 > |z|^2 \geq 400.0$ | *6 convergents (Laplace continued fractions)* |
|  | $400.0 > |z|^2 \geq 127.0$ | *9 terms of the series in (4)* |
| V | $127.0 > |z|^2$ & $y^2 \leq 10^{-3}$ | *Daw715(x) and a few terms of Taylor series* |
| VI | *Otherwise* | *Algorithm 916 (loop used for $x < (-\ln R_{min})^{1/2}$)* |

## 3. Speed Benchmark

The present improvements in calculating the complex probability function have been implemented in a *Fortran* elemental subroutine which can be run using single or double precision arithmetic. An efficiency comparison is performed between the present CPF code and the *Fortran* version of Algorithm 916, which included efficiency improvements compared to the original version of the algorithm. The comparison has been performed for the same four datasets used for performance comparison in Zaghloul 2016 and Zaghloul 2017[2] each comprising a total of 2,840,071 input points in the complex domain. Each case has been evaluated for 100 consecutive times and the average run time is used to report the time per evaluation. The results of the comparison, using double precision (*Fortran-d)* and single precision (*Fortran-s*) are given in Table 4. As can be seen from the table, the present algorithm is consistently faster than the latest version of Algorithm 916 when run using double or single precision arithmetic. Depending on the case studied, targeted accuracy, and the precision used, efficiency improvements vary and can go up to more than a factor of five. For case 3 and case 4, one can calculate the function, using the present algorithm, to 13 significant digits, in a shorter time than used to calculate it to 5 significant digits with Algorithm 916.

---

[2] The *y* values are 71 values that are logarithmically equally spaced between $10^{-5}$ and $10^5$ for Case 1; $10^{-20}$ and $10^4$ for Case 2; $10^{-5}$ and $10^5$ for Case 3; and $10^{-20}$ and 6 for Case 4. The *x* values are 40001 values linearly equally spaced between -500 and 500 for Case 1; -200 and 200 for Case 2; -10 and 10 for Case 3; and *x randomly generated with $|z|^2 \leq 36$ for Case 4.*

A similar comparison is held with the free/open source *"Faddeeva Package"* developed in *C++* at the *Massachusetts Institute of Technology* (*MIT*) by S. Johnson [Johnson 2012]**.** The package uses a combination of Algorithm 916 [Zaghloul and Ali 2011] for relatively small $|z|$ and Algorithm 680 [Poppe and Wijers. 1990] or more precisely continued fraction approximations for large values of $|z|$. Although we use MIT-C++ package herein for speed benchmarking, it has to be noted that the package does not satisfy its claimed accuracy where the accuracy deteriorates down to 6, 7 or 8 significant figures while requesting 13 significant digits (through setting the required relative error to 0.0). Just for example, we refer here to a few points where the accuracy of the MIT-C++ package in calculating the Faddeyeva function is declined as mentioned. Consider for example the points ($z=\pm6.0+0.158489319246111i$, $z=\pm6.0+0.138949549437314i$, $z=0.0+7.19685673001151i$, $z=0.0+8.20891415963826i$) and nearby points. Although it seems to be unfair efficiency benchmark, we hold it as it gives more credibility to the efficiency of the present code.

The version of the package "*Faddeeva.cc*" used in this comparison was last modified on May 12 2015. The results of this comparison are shown in Table 5. For all cases, the present algorithm is consistently and considerably faster for all targeted accuracies. Depending on the case studied and the targeted accuracy, the present algorithm can be up to a factor of five faster than the *MIT-C++* package. Even though, it is recognized that the run time for the *MIT-C++* package does not correlate with the targeted accuracy, in a clear contradiction with the logical expectations. For the case of 13 significant digits, Table 5 reveals that the present algorithm can be faster by a factor greater than 2 depending on the case under consideration. This factor increases to 5 for the case of 4 significant digits.

**Table 4.** Speed comparison between *Fortran* implementations of the present routine and the latest version of Algorithm 916. The values have been generated using Intel Visual Fortran 64 Compiler Professional for applications running on Intel(R) 64, Version 11.1.

|  | Average time per evaluation (ns) | | | | | | | | | | | |
|---|---|---|---|---|---|---|---|---|---|---|---|---|
|  | Case 1 | | | Case 2 | | | Case 3 | | | Case 4 | | |
| Algorithm | 916-V2 | present | ratio | 916-V2 | present | Ratio | 916-V2 | present | ratio | 916-V2 | present | ratio |
| **Double Precision** | | | | | | | | | | | | |
| *Faddeyeva (z, 13)* | 39.6 | 15.2 | 0.38 | 105.5 | 18.6 | 0.18 | 158.4 | 74.8 | 0.47 | 222.5 | 77.7 | *0.35* |
| *Faddeyeva (z, 12)* | 30.0 | 14.5 | 0.48 | 76.5 | 17.4 | 0.23 | 147.8 | 72.2 | 0.49 | 215.0 | 73.6 | *0. 34* |
| *Faddeyeva (z, 11)* | 24.0 | 14.1 | 0.59 | 57.7 | 16.8 | 0.29 | 140.3 | 69.2 | 0.49 | 208.1 | 72.3 | *0.35* |
| *Faddeyeva (z, 10)* | 20.7 | 13.5 | 0.65 | 47.0 | 16.6 | 0.35 | 137.5 | 68.9 | 0.50 | 208.1 | 71.8 | *0.35* |
| *Faddeyeva (z, 9)* | 17.9 | 11.5 | 0.64 | 38.2 | 16.3 | 0.43 | 130.3 | 66.3 | 0.51 | 200.5 | 70.3 | *0.35* |
| *Faddeyeva (z, 8)* | 16.3 | 10.6 | 0.65 | 33.4 | 14.4 | 0.43 | 127.5 | 65.7 | 0.52 | 200.4 | 69.2 | *0.35* |
| *Faddeyeva (z, 7)* | 13.6 | 10.0 | 0.74 | 24.5 | 12.8 | 0.52 | 118.8 | 60.1 | 0.51 | 193.0 | 65.9 | *0.34* |
| *Faddeyeva (z, 6)* | 12.7 | 9.5 | 0.75 | 19.8 | 12.1 | 0.61 | 115.5 | 48.5 | 0.42 | 193.6 | 60.6 | *0.31* |
| *Faddeyeva (z, 5)* | 12.6 | 9.5 | 0.75 | 16.9 | 10.7 | 0.63 | 109.2 | 29.7 | 0.27 | 186.2 | 54.5 | *0.29* |
| *Faddeyeva (z, 4)* | 10.9 | 8.5 | 0.78 | 11.7 | 9.5 | 0.81 | 21.3 | 20.6 | 0.97 | 50.8 | 43.7 | *0.86* |
| | | | | | | | | | | | | |
| **Single Precision** | | | | | | | | | | | | |
| *Faddeyeva (z, 6)* | 7.9 | 7.7 | 0.97 | 11.4 | 9.2 | 0.81 | 65.7 | 38.5 | 0.59 | 110.6 | 49.5 | *0.45* |
| *Faddeyeva (z, 5)* | 7.6 | 6.5 | 0.86 | 10.6 | 8.1 | 0.76 | 59.9 | 23.7 | 0.40 | 102.3 | 39.8 | *0.39* |
| *Faddeyeva (z, 4)* | 6.5 | 6.0 | 0.92 | 8.4 | 7.3 | 0.87 | 17.6 | 16.4 | 0.93 | 43.8 | 35.7 | *0.82* |

**Table 5.** Speed comparison between the *Fortran* implementation of the present routine and the *MIT-C++ "Faddeeva Package"*. The values have been generated using Intel Visual Fortran 64 Compiler Professional for applications running on Intel(R) 64, Version 11.1 and Microsoft Visual C++ 2015 x64 x86 Build Tools, Optimizing Compiler Version 19.00.24210 for x86.

|  | Average time per evaluation (ns) | | | | | | | | | | | |
|---|---|---|---|---|---|---|---|---|---|---|---|---|
|  | Case 1 | | | Case 2 | | | Case 3 | | | Case 4 | | |
| Algorithm | MIT C++ | present DP | ratio | MIT C++ | present DP | ratio | MIT C++ | present DP | ratio | MIT C++ | present DP | Ratio |
| *relerr= <1e-13* | 30.9 | 15.2 | 0.49 | 49.0 | 18.6 | 0.38 | 97.9 | 74.8 | 0.76 | 104.8 | 77.7 | *0.74* |
| *relerr= <1e-12* | 30.1 | 14.5 | 0.48 | 48.0 | 17.4 | 0.36 | 98.8 | 72.2 | 0.73 | 105.2 | 73.6 | *0.70* |
| *relerr= <1e-11* | 30.3 | 14.1 | 0.47 | 48.2 | 16.8 | 0.35 | 96.7 | 69.2 | 0.72 | 103.5 | 72.3 | *0.70* |
| *relerr= <1e-10* | 30.2 | 13.5 | 0.45 | 47.5 | 16.6 | 0.35 | 98.1 | 68.9 | 0.70 | 104.5 | 71.8 | *0.69* |
| *relerr= <1e-09* | 30.1 | 11.5 | 0.38 | 47.4 | 16.3 | 0.34 | 97.1 | 66.3 | 0.68 | 103.5 | 70.3 | *0.68* |
| *relerr= <1e-08* | 30.3 | 10.6 | 0.35 | 47.6 | 14.4 | 0.30 | 97.8 | 65.7 | 0.67 | 103.4 | 69.2 | *0.67* |
| *relerr= <1e-07* | 30.1 | 10.0 | 0.33 | 47.6 | 12.8 | 0.27 | 96.8 | 60.1 | 0.62 | 103.6 | 65.9 | *0.64* |
| *relerr= <1e-06* | 30.4 | 9.5 | 0.31 | 47.6 | 12.1 | 0.25 | 97.8 | 48.5 | 0.50 | 103.9 | 60.6 | *0.58* |
| *relerr= <1e-05* | 30.1 | 9.5 | 0.32 | 47.5 | 10.7 | 0.22 | 97.2 | 29.7 | 0.31 | 103.4 | 54.5 | *0.53* |
| *relerr= <1e-04* | 30.2 | 8.5 | 0.28 | 47.7 | 9.5 | 0.20 | 98.1 | 20.6 | 0.21 | 103.9 | 43.7 | *0.42* |

## 4. Computation of Related Special Functions with Complex Arguments

The complex probability function is related to the complementary error function *erfc(z)*, the error function *erf(z)*, the imaginary error function *erfi(z)*, and the scaled complementary error function *erfcx(z)* as given in Eq. (1), and one can simply calculate these functions using the following direct relations

$$erfc(z) = e^{-z^2} w(iz)$$
$$erf(z) = 1 - erfc(z) = 1 - e^{-z^2} w(iz)$$
$$erfi(z) = i\, erf(-iz) = i(1 - e^{z^2} w(z))$$
$$erfcx(z) = e^{z^2} erfc(z) = w(iz)$$

(7)

Both of *erf(z)* and *erfi(z)* have converging series expansions near *z=0* [Abramowitz and Stegun, 1964] which are used herein for the computation of the functions for $|z| \leq 1$. It is easy to show that

$$erfc(iy) = 1 - i\, e^{y^2} \frac{2}{\sqrt{\pi}} Daw(y)$$
$$erf(iy) = i\, e^{y^2} \frac{2}{\sqrt{\pi}} Daw(y)$$
$$erfi(x) = e^{y^2} \frac{2}{\sqrt{\pi}} Daw(x)$$

(8)

Similarly, the Dawson's function or Dawson's integral, *Daw(z)*, and the plasma dispersion *Zeta* function, *Z(z)*, can be easily found from the complex probability function using the relations

$$Daw(z) = \frac{\sqrt{\pi}}{2} e^{-z^2} erfi(z) = \frac{i\sqrt{\pi}}{2} (e^{-z^2} - w(z))$$
$$Z(z) = i\sqrt{\pi}\, w(z)$$

(9)

However, *Daw(x)* from Algorithm 715 and a few terms of the Taylor series have been used to compute *Daw(z)* directly for small values of *y* and non-negligible values of *exp(-z²)*.

Perhaps the interesting and relatively involved part is the computation of the sine and cosine Fresnel integrals, which are defined in a variety of equivalent forms in the mathematical literature (see for example Abramowitz and Stegun 1964 and Zwillinger 2003). For example, the functions *S(z)* and *C(z)* are defined by the integrals

$$S(z) = \int_0^z \sin(\frac{\pi}{2} t^2)\, dt$$
$$C(z) = \int_0^z \cos(\frac{\pi}{2} t^2)\, dt$$

(10)

The functions are also defined in the forms of the pairs $S_1(z)$, $C_1(z)$ and $S_2(z)$, $C_2(z)$ where

$$S_1(z) = \sqrt{\frac{\pi}{2}} \int_0^z \sin(t^2)\,dt$$

$$C_1(z) = \sqrt{\frac{\pi}{2}} \int_0^z \cos(t^2)\,dt \tag{11}$$

$$S_2(z) = \frac{1}{\sqrt{2\pi}} \int_0^z \frac{\sin(t)}{\sqrt{t}}\,dt = \frac{1}{2} \int_0^z J_{1/2}(t)\,dt$$

$$C_2(z) = \frac{1}{\sqrt{2\pi}} \int_0^z \frac{\cos(t)}{\sqrt{t}}\,dt = \frac{1}{2} \int_0^z J_{-1/2}(t)\,dt \tag{12}$$

where $J_{1/2}(t)$ and $J_{-1/2}(t)$ are the ordinary Bessel functions of the first kind of orders ½ and -½ respectively. The three pairs of functions in (10), (11) and (12) are related to each other by

$$S(z) = S_1\left(\sqrt{\tfrac{\pi}{2}}\,z\right) = S_2\left(\tfrac{\pi}{2}z^2\right)$$

$$C(z) = C_1\left(\sqrt{\tfrac{\pi}{2}}\,z\right) = C_2\left(\tfrac{\pi}{2}z^2\right) \tag{13}$$

From the computational point of view, all forms of the Fresnel integrals can be computed using the same algorithm with a slight change in the argument as explained in Eq. (13).

Considering the form of the integrals given in (10), it is easy to show that they can be also expressed in terms of the complex probability function as given below;

$$S(z) = \frac{-1-i}{4}\left[\left(1 - e^{(\tfrac{1-i}{2}\sqrt{\pi}z)^2} w\left(\frac{1-i}{2}\sqrt{\pi}z\right)\right) + i\left(1 - e^{(\tfrac{1+i}{2}\sqrt{\pi}z)^2} w\left(\frac{1+i}{2}\sqrt{\pi}z\right)\right)\right] \tag{14}$$

and

$$C(z) = \frac{-1+i}{4}\left[\left(1 - e^{(\tfrac{1-i}{2}\sqrt{\pi}z)^2} w\left(\tfrac{1-i}{2}\sqrt{\pi}z\right)\right) - i\left(1 - e^{(\tfrac{1+i}{2}\sqrt{\pi}z)^2} w\left(\tfrac{1+i}{2}\sqrt{\pi}z\right)\right)\right] \tag{15}$$

Fresnel integrals are analytical functions defined for all complex values of $z$, over the whole complex plane. They are odd functions with the following symmetry relations

$$\begin{aligned} S(-z) &= -S(z) & S(\bar{z}) &= \overline{S(z)} \\ C(-z) &= -C(z) & C(\bar{z}) &= \overline{C(z)} \end{aligned} \tag{16}$$

In addition, the integrals have simple limiting values at $z$=zero and $z\to\infty$ where

$$S(0) = 0, \quad \lim_{z \to \pm\infty} S(z) = \pm\frac{1}{2},$$
$$C(0) = 0, \quad \lim_{z \to \pm\infty} C(z) = \pm\frac{1}{2}, \tag{17}$$

It is also worth mentioning that the Fresnel integrals $S(z)$ and $C(z)$ have the following simple converging series representations near $z \to 0$,

$$S(z) = z^3 \sum_{k=0}^{\infty} \frac{2^{-2k-1} \pi^{2k+1} (-z^4)^k}{(4k+3)(2k+1)!}$$
$$C(z) = z \sum_{k=0}^{\infty} \frac{2^{-2k} \pi^{2k} (-z^4)^k}{(4k+1)(2k)!} \tag{18}$$

Using the expressions (14) and (15) with some mathematical manipulation, one can show that

$$S(x) = -\tfrac{1}{2} + \tfrac{1}{2}\mathrm{Re}(w(\tfrac{1-i}{2}\sqrt{\pi}\, x)) \times \left(\cos(\tfrac{\pi}{2}x^2) + \sin(\tfrac{\pi}{2}x^2)\right)$$
$$\quad + \tfrac{1}{2}\mathrm{Im}(w(\tfrac{1-i}{2}\sqrt{\pi}\, x)) \times \left(\sin(\tfrac{\pi}{2}x^2) - \cos(\tfrac{\pi}{2}x^2)\right)$$
with
$$S(iy) = -i\,S(y) \tag{19}$$

Similarly,

$$C(x) = -\tfrac{1}{2} + \tfrac{1}{2}\mathrm{Re}(w(\tfrac{1-i}{2}\sqrt{\pi}\, x)) \times \left(\cos(\tfrac{\pi}{2}x^2) - \sin(\tfrac{\pi}{2}x^2)\right)$$
$$\quad + \tfrac{1}{2}\mathrm{Im}(w(\tfrac{1-i}{2}\sqrt{\pi}\, x)) \times \left(\cos(\tfrac{\pi}{2}x^2) + \cos(\tfrac{\pi}{2}x^2)\right)$$
with
$$C(i\,y) = i\,C(y) \tag{20}$$

The formulation of the functions as direct expressions in the complex probability function as in Eqs. (14,15) and in Eqs. (19,20) is not familiar in the literature although equivalent alternatives can be found.

The expansions in (18) are used herein for the computation of the functions for $|z|\leq 1$, while the expressions in (19) and (20) are used together with a few terms from the Taylor series to calculate the functions for $|z|>1$ with $y \leq 10^{-4}$ and/or $x \leq 10^{-4}$ while (14) and (15) are used for the rest of the domain. A special care has been devoted to the calculation of the cosine and sine functions in (19) and (20) when the argument is an integer value of $\pi/2$ to avoid deterioration of the accuracy due to rounding error.

## 5. Accuracy and Representative Results for Related Special Functions

In the appendix section, Tables 6, 7, 8 and 9 show representative results for $erfc(z)$, $erf(z)$, $erfi(z)$ and $Daw(z)$ as calculated from the present algorithm compared to values calculated from the *Matlab Special Functions of Applied Mathematics (SFAM) in the symbolic toolbox* and from the *Mathematica* and *Maple* software packages. Computations of these functions embody computations of either $exp(z^2)$ or $exp(-z^2)$ in addition to the calculation of $w(z)$ or $w(u(z))$ where $u(z)$ is a complex argument that depends on $z$. Numeric overflow and underflow occur during computations of the term $exp(z^2)$ for $x^2-y^2>log(R_{max})$ and $y^2-x^2>log(R_{max})$, in order, where $R_{max}$ is the largest floating point number in the computational platform. Similarly, numeric overflow and underflow occur during computations of the term $exp(-z^2)$ for $y^2-x^2>log(R_{max})$ and $x^2-y^2>log(R_{max})$, respectively. On the other hand, computations of $w(u(z))$ can fail as a result of unavoidable overflow for arguments of large magnitudes with negative imaginary part. Whenever such a failure due to a possible overflow could occur, the code returns (NaN+NaN i) without experiencing an overflow.

Investigation of Tables 6-9 shows that results from the present work and those from the *Matlab SFAM* and the *Maple* software package are in very good agreement (up to 13 significant digits). However, there are some points, where results from *Mathematica* are not at that level of agreement (see for example results for *erfc* and *erf* for $z=5.9\times10^{-10}+15.0i$, and results for *erfi* and *Daw* for $z=6.3+1\times10^{-10}i$) while results from *Mathematica* for some other points seem to be anomalous (see for example results for $z=6.3\times10^{-10}+25.0i$).

Accuracy comparison for the calculation of Fresnel function $S(z)$ can be found in Tables 10. A very good agreement (equal to or better than 10 significant digits) between the present results and those calculated using *Matlab* and *Maple* is recognizable with a single exception for the point $z=26.0+1\times10^{-2}i$ where the agreement for the imaginary part deteriorates to only 8 significant digits with *Matlab*, 9 significant digits with *Maple*, and 10 significant digits with *Mathematica*. Again, there are many points where results from *Mathematica* show deteriorated accuracy or anomalous results, which is very clear for $z=6.3\times10^{-10}+25.0i$.

For results of $C(z)$, as shown in Table 11, the agreement with *Matlab* calculations is reasonable (equal to or better than 9 significant digits) except for the real part of the function for the point $z=5.9\times10^{-10}+15.0i$ where there is no such agreement for the real part. For this point there is also a disagreement between the results of *the present work* and those of *Mathematica* and *Maple*. However, repeating calculations with *Maple* using 32 digits, it was found that agreement up to 10 significant digits for the real part of this point is obtained. A similar discrepancy between results for the imaginary part from the present work and those from *Maple* for the point $z=23.0+1.0\times10^{-5}i$ is resolved upon repeating *Maple* calculations to 32 digits. This undoubtedly gives credibility to the accuracy of the present package. The calculation of $erfcx(z)$ and $Z(z)$ is effectively computation of the complex probability function.

## 6. CONCLUSIONS

An algorithm for efficient multi-accuracy computations of a set of special functions with complex argument is presented. The present algorithm showed superior efficiency and accuracy compared to standard algorithms in the literature. Using results from the present algorithm, concerns about in-accuracies in widely used free and commercial software packages have been pointed out and highlighted.

# Appendix

```
*******************************************
Table 6.  Comparison between values of erfc(z) calculated
from the present package and those calculated from
 Matlab (9.2.0.538062(R2017a)), Mathematica9 and Maple13
*******************************************
     x         y          Re(Matlab)               Im(Matlab)              Re(Mathematica)          Im(Mathematica)           Re(Maple)               Im(Maple)               Re(present)             Im(present)
 0.630D-09  0.100D-09  0.9999999992891211D+00 -0.1128379167095513D-09  0.9999999992891211D+00 -0.1128379167095513D-09  0.9999999992891211D+00 -0.1128379167095513D-09  0.9999999992891211D+00 -0.1128379167095513D-09
 0.630D-08  0.100D-09  0.9999999928912112D+00 -0.1128379167095513D-09  0.9999999928912112D+00 -0.1128379167095513D-09  0.9999999928912112D+00 -0.1128379167095513D-09  0.9999999928912112D+00 -0.1128379167095512D-09
 0.230D-06  0.100D-09  0.9999997404727916D+00 -0.1128379167095453D-09  0.9999997404727916D+00 -0.1128379167095453D-09  0.9999997404727916D+00 -0.1128379167095453D-09  0.9999997404727916D+00 -0.1128379167095453D-09
 0.630D-01  0.100D-09  0.9290060498694527D+00 -0.1123909506091105D-09  0.9290060498694527D+00 -0.1123909506091105D-09  0.9290060498694527D+00 -0.1123909506091105D-09  0.9290060498694528D+00 -0.1123909506091106D-09
 0.630D+00  0.100D-08  0.3729535566618043D+00 -0.7587235764175012D-09  0.3729535566618043D+00 -0.7587235764175011D-09  0.3729535566618043D+00 -0.7587235764175011D-09  0.3729535566618042D+00 -0.7587235764175011D-09
 0.630D+01  0.100D-07  0.5124221687395676D-18 -0.6535925189178454D-25  0.5124221687395672D-18 -0.6535925189178451D-25  0.5124221687395663D-18 -0.6535925189178438D-25  0.5124221687395672D-18 -0.6535925189178449D-25
 0.230D+02  0.100D-04  0.4441265477758837D-231 -0.2044909616141931D-234  0.4441265477758918D-231 -0.2044909616141968D-234  0.4441265477758837D-231 -0.2044909616141931D-234  0.4441265477759033D-231 -0.2044909616142021D-234
 0.510D-09  0.250D+00  0.9999999993874118D+00 -0.2880836197949720D+00  0.9999999993874118D+00 -0.2880836197949721D+00  0.9999999993874118D+00 -0.2880836197949720D+00  0.9999999993874117D+00 -0.2880836197949720D+00
 0.590D-09  0.150D+02 -0.3463901223474738D+89 -0.1961384563867380D+97 -0.3463901219342600D+89 -0.1961384563867408D+97 -0.3463901223474738D+89 -0.1961384563867380D+97 -0.3463901223474737D+89 -0.1961384563867380D+97
 0.630D-09  0.250D+02 -0.1931286916174148D+263 -0.6135986249821948D+270 -0.2721857454564708D+285 -0.8314637164730795D+292 -0.1931286916174148D+263 -0.6135986249821948D+270 -0.1931286916174148D+263 -0.6135986249821947D+270
 0.630D-07  0.100D-09  0.9999999289121125D+00 -0.1128379167095508D-09  0.9999999289121124D+00 -0.1128379167095508D-09  0.9999999289121125D+00 -0.1128379167095508D-09  0.9999999289121124D+00 -0.1128379167095508D-09
 0.630D-07  0.100D-06  0.9999999289121125D+00 -0.1128379167095512D-06  0.9999999289121124D+00 -0.1128379167095512D-06  0.9999999289121125D+00 -0.1128379167095512D-06  0.9999999289121124D+00 -0.1128379167095512D-06
 0.630D+01  0.100D-09  0.5124221687395717D-18 -0.6535925189178471D-27  0.5124221687395704D-18 -0.6535925189178453D-27  0.5124221687395715D-18 -0.6535925189178466D-27  0.5124221687395715D-18 -0.6535925189178468D-27
 0.630D+01  0.100D-07  0.5124221687395676D-18 -0.6535925189178454D-25  0.5124221687395672D-18 -0.6535925189178451D-25  0.5124221687395663D-18 -0.6535925189178438D-25  0.5124221687395672D-18 -0.6535925189178449D-25
 0.630D+01  0.100D-03  0.5124217569763372D-18 -0.6535923481559211D-21  0.5124217569763389D-18 -0.6535923481559230D-21  0.5124217569763359D-18 -0.6535923481559194D-21  0.5124217569763374D-18 -0.6535923481559212D-21
 0.630D+01  0.100D+01  0.1352712477970645D-17 -0.2565734713428141D-18  0.1352712477970645D-17 -0.2565734713428137D-18  0.1352712477970641D-17 -0.2565734713428140D-18  0.1352712477970644D-17 -0.2565734713428141D-18
 0.630D+01  0.250D+02 -0.1927240394121577D+253 -0.2856604986670404D+253 -0.1927240394121470D+253 -0.2856604986670310D+253 -0.1927240394121601D+253 -0.2856604986670377D+253 -0.1927240394121496D+253 -0.2856604986670222D+253
 0.630D+01  0.260D+02 -0.2781086908548881D+275 -0.3764149461796603D+275 -0.2781086908548657D+275 -0.3764149461796397D+275 -0.2781086908548913D+275 -0.3764149461796564D+275 -0.2781086908548633D+275 -0.3764149461796453D+275
 0.630D+01  0.120D+02  0.6846033501217893D+43 -0.8314590451931710D+44  0.6846033501217729D+43 -0.8314590451931671D+44  0.6846033501217479D+43 -0.8314590451931692D+44  0.6846033501218439D+43 -0.8314590451931703D+44
 0.130D+02  0.150D+02  0.1602309039170805D+23 -0.5726318806999969D+23  0.1602309039170810D+23 -0.5726318806999964D+23  0.1602309039170805D+23 -0.5726318806999969D+23  0.1602309039170806D+23 -0.5726318806999968D+23
 0.130D+02  0.200D+02 -0.2702312808779561D+172  0.2674523742514406D+171 -0.2702312808779589D+172  0.2674523742514516D+171 -0.2702312808779561D+172  0.2674523742514406D+171 -0.2702312808779561D+172  0.2674523742514462D+171
 0.230D+01  0.100D+02 -0.7305679310276024D+40  0.1623858803164295D+40 -0.7305679310276026D+40  0.1623858803164313D+40 -0.7305679310276012D+40  0.1623858803164322D+40 -0.7305679310276067D+40  0.1623858803164334D+40
 0.260D+02  0.100D-01  0.4914036960738840D-295 -0.2816096460526261D-295  0.4914036960738912D-295 -0.2816096460526302D-295  0.4914036960738840D-295 -0.2816096460526261D-295  0.4914036960526190D-295
 0.730D+01  0.100D-03  0.5502567832011327D-24 -0.8107780186742224D-27  0.5502567832011307D-24 -0.8107780186742195D-27  0.5502567832011311D-24 -0.8107780186742200D-27  0.5502567832011319D-24 -0.8107780186742213D-27
 0.630D+01  0.100D-06  0.5124221687391599D-18 -0.6535925189176764D-24  0.5124221687391584D-18 -0.6535925189176741D-24  0.5124221687391586D-18 -0.6535925189176747D-24  0.5124221687391584D-18 -0.6535925189176744D-24
 0.630D+01  0.100D-05  0.5124221686983954D-18 -0.6535925189007709D-23  0.5124221686983971D-18 -0.6535925189007729D-23  0.5124221686983941D-18 -0.6535925189007693D-23  0.5124221686983963D-18 -0.6535925189007721D-23
 0.430D+00  0.100D-05  0.5431133054500564D+00 -0.9378945443478894D-06  0.5431133054500564D+00 -0.9378945443478894D-06  0.5431133054500564D+00 -0.9378945443478894D-06  0.5431133054500565D+00 -0.9378945443478894D-06
 0.130D+02  0.100D-01  0.1680915860284974D-74 -0.4485365503481142D-75  0.1680915860284967D-74 -0.4485365503481123D-75  0.1680915860284974D-74 -0.4485365503481142D-75  0.1680915860284980D-74 -0.4485365503481156D-75
 0.430D+01  0.100D+02  0.9055021941027503D+34  0.9379608607850867D+34  0.9055021941027523D+34  0.9379608607850810D+34  0.9055021941027524D+34  0.9379608607850814D+34  0.9055021941027573D+34  0.9379608607850862D+34
 0.630D-01  0.120D+02 -0.1620124184557827D+62 -0.1039653977037593D+61 -0.1620124184557789D+62 -0.1039653977037569D+61 -0.1620124184557827D+62 -0.1039653977037593D+61 -0.1620124184557807D+62 -0.1039653977037580D+61
 0.630D-01  0.150D+02 -0.1857482893618806D+97  0.6052285795414974D+96 -0.1857482893618797D+97  0.6052285795414943D+96 -0.1857482893618806D+97  0.6052285795414974D+96 -0.1857482893618783D+97  0.6052285795414903D+96
 0.630D-01  0.200D+02 -0.8591804934512793D+172  0.1191477261295850D+173 -0.8591804934512993D+172  0.1191477261295879D+173 -0.8591804934512793D+172  0.1191477261295850D+173 -0.8591804934512906D+172  0.1191477261295866D+173
 0.130D+02  0.160D+02 -0.9749075032820455D+36 -0.1347556133349446D+37 -0.9749075032820759D+36 -0.1347556133349415D+37 -0.9749075032820455D+36 -0.1347556133349446D+37 -0.9749075032820488D+36 -0.1347556133349445D+37
 0.100D+01  0.100D-01  0.1572576956087019D+00 -0.4150936598121530D-02  0.1572576956087018D+00 -0.4150936598121530D-02  0.1572576956087019D+00 -0.4150936598121530D-02  0.1572576956087018D+00 -0.4150936598121518D-02
 0.550D+01  0.100D-03  0.7357843395151008D-14 -0.8223314414305329D-17  0.7357843395151023D-14 -0.8223314414305349D-17  0.7357843395151008D-14 -0.8223314414305329D-17  0.7357843395151009D-14 -0.8223314414305334D-17
 0.190D+01  0.100D-01  0.1252479183851085D-157 -0.4619188896571069D-158  0.1252479183851073D-157 -0.4619188896570986D-158  0.1252479183851085D-157 -0.4619188896571069D-158  0.1252479183851085D-157 -0.4619188896571069D-158
 0.100D+01  0.280D+02  Infinity               -Infinity                Infinity               -Infinity                Infinity               -Infinity                Infinity                NaN
 0.150D+02  0.150D-01 -0.9109691190248829D-03  0.2658046409880405D-01 -0.9109691190241991D-03  0.2658046409880407D-01 -0.9109691190248829D-03  0.2658046409880405D-01 -0.9109691190248848D-03  0.2658046409880409D-01
 0.260D+02  0.000D+00  0.5663192408856143D-295  0.0000000000000000D+00  0.5663192408856014D-295  0.0000000000000000D+00  0.5663192408856143D-295  0.0000000000000000D+00  0.5663192408856143D-295  0.0000000000000000D+00
 0.000D+00  0.266D+02  0.1000000000000000D+01 -0.4132896053051959D+306  0.2530668961289347D+290 -0.4132896053051863D+306  0.1000000000000000D+01 -0.4132896053051739D+306  0.1000000000000000D+01 -0.4132896053051984D+306
```

```
******************************************
******************************************
Table 7.  Comparison between values of erf(z) calculated
from the present package and those calculated from
 Matlab (9.2.0.538062(R2017a)), Mathematica9 and Maple13
******************************************
     x           y        Re(Matlab)              Im(Matlab)             Re(Mathematica)         Im(Mathematica)          Re(Maple)              Im(Maple)              Re(present)            Im(present)
 0.630D-09   0.100D-09   0.7108788752701729D-09  0.1128379167095513D-09  0.7108788752701730D-09  0.1128379167095513D-09  0.7108788752701729D-09  0.1128379167095513D-09  0.7108788752701729D-09  0.1128379167095513D-09
 0.630D-08   0.100D-09   0.7108788752701729D-08  0.1128379167095513D-09  0.7108788752701730D-08  0.1128379167095513D-09  0.7108788752701729D-08  0.1128379167095513D-09  0.7108788752701730D-08  0.1128379167095513D-09
 0.230D-06   0.100D-09   0.2595272084319633D-06  0.1128379167095453D-09  0.2595272084319633D-06  0.1128379167095453D-09  0.2595272084319633D-06  0.1128379167095453D-09  0.2595272084319633D-06  0.1128379167095453D-09
 0.630D-01   0.100D-09   0.7099395013054732D-01  0.1123909506091105D-09  0.7099395013054732D-01  0.1123909506091105D-09  0.7099395013054732D-01  0.1123909506091105D-09  0.7099395013054732D-01  0.1123909506091105D-09
 0.630D+00   0.100D-08   0.6270464433381957D+00  0.7587235764175012D-09  0.6270464433381957D+00  0.7587235764175011D-09  0.6270464433381957D+00  0.7587235764175011D-09  0.6270464433381948D+00  0.7587235764174664D-09
 0.630D+01   0.100D-07   0.1000000000000000D+01  0.6535925189178454D-25  0.1000000000000000D+01  0.6535925189178451D-25  0.1000000000000000D+01  0.6535925189178438D-25  0.1000000000000000D+01  0.6535925189178449D-25
 0.230D+02   0.100D-04   0.1000000000000000D+01  0.2044909616141931-234  0.1000000000000000D+01  0.2044909616141968-234  0.1000000000000000D+01  0.2044909616141931-234  0.1000000000000000D+01  0.2044909616142021-234
 0.510D-09   0.250D+00   0.6125882191750765D-09  0.2880836197949720D+00  0.6125882191750765D-09  0.2880836197949721D+00  0.6125882191750765D-09  0.2880836197949720D+00  0.6125882191750764D-09  0.2880836197949719D+00
 0.590D-09   0.150D+02   0.3463901223474738D+89  0.1961384563867380D+97  0.3463901219342600D+89  0.1961384563867408D+97  0.3463901223474738D+89  0.1961384563867380D+97  0.3463901223474737D+89  0.1961384563867380D+97
 0.630D-09   0.250D+02   0.1931286916174148+263  0.6135986249821948+270  0.2721857454564708+285  0.8314637164730795+292  0.1931286916174148+263  0.6135986249821948+270  0.1931286916174148+263  0.6135986249821947+270
 0.630D-07   0.100D-07   0.7108788752701718D-07  0.1128379167095508D-09  0.7108788752701720D-07  0.1128379167095508D-09  0.7108788752701720D-07  0.1128379167095508D-09  0.7108788752701720D-07  0.1128379167095508D-09
 0.630D-07   0.100D-06   0.7108788752701790D-07  0.1128379167095512D-06  0.7108788752701793D-07  0.1128379167095512D-06  0.7108788752701791D-07  0.1128379167095512D-06  0.7108788752701791D-07  0.1128379167095512D-06
 0.630D+01   0.100D-07   0.1000000000000000D+01  0.6535925189178471D-27  0.1000000000000000D+01  0.6535925189178454D-25  0.1000000000000000D+01  0.6535925189178451D-25  0.1000000000000000D+01  0.6535925189178438D-25  0.1000000000000000D+01  0.6535925189178449D-25
 0.630D+01   0.100D-03   0.1000000000000000D+01  0.6535923481559211D-21  0.1000000000000000D+01  0.6535923481559230D-21  0.1000000000000000D+01  0.6535923481559194D-21  0.1000000000000000D+01  0.6535923481559212D-21
 0.630D+01   0.100D+03   0.1000000000000000D+01  0.2565734713428141D-18  0.1000000000000000D+01  0.2565734713428137D-18  0.1000000000000000D+01  0.2565734713428140D-18  0.1000000000000000D+01  0.2565734713428141D-18
 0.630D+01   0.250D+02   0.1927240394121577+253  0.2856604986670404+253  0.1927240394121470+253  0.2856604986670310+253  0.1927240394121601+253  0.2856604986670377+253  0.1927240394121496+253  0.2856604986670222+253
 0.630D+01   0.260D+02   0.2781086908548881+275  0.3764149461796603+275  0.2781086908548657+275  0.3764149461796397+275  0.2781086908548913+275  0.3764149461796564+275  0.2781086908548633+275  0.3764149461796453+275
 0.630D+01   0.120D+02  -0.6846033501217893D+43  0.8314590451931710D+44 -0.6846033501217729D+43  0.8314590451931670D+44 -0.6846033501217479D+43  0.8314590451931692D+44 -0.6846033501218439D+43  0.8314590451931703D+44
 0.130D+02   0.150D+02 -0.1602309039170805D+23  0.5726318806999969D+23 -0.1602309039170810D+23  0.5726318806999964D+23 -0.1602309039170805D+23  0.5726318806999969D+23 -0.1602309039170806D+23  0.5726318806999968D+23
 0.130D+02   0.200D+02  0.2702312808779561+172 -0.2674523742514406+171  0.2702312808779589+172 -0.2674523742514516+171  0.2702312808779561+172 -0.2674523742514406+171  0.2702312808779561+172 -0.2674523742514462+171
 0.230D+01   0.100D+02  0.7305679310276024D+40 -0.1623858803164295D+40  0.7305679310276026D+40 -0.1623858803164313D+40  0.7305679310276012D+40 -0.1623858803164322D+40  0.7305679310276067D+40 -0.1623858803164334D+40
 0.260D+02   0.100D-01   0.1000000000000000D+01  0.2816096460526261-295  0.1000000000000000D+01  0.0000000000000000D+00  0.1000000000000000D+01  0.2816096460526261-295  0.1000000000000000D+01  0.2816096460526190-295
 0.730D+01   0.100D-03   0.1000000000000000D+01  0.8107780186742224D-27  0.1000000000000000D+01  0.8107780186742195D-27  0.1000000000000000D+01  0.8107780186742200D-27  0.1000000000000000D+01  0.8107780186742213D-27
 0.630D+01   0.100D-06   0.1000000000000000D+01  0.6535925189176764D-24  0.1000000000000000D+01  0.6535925189176741D-24  0.1000000000000000D+01  0.6535925189176747D-24  0.1000000000000000D+01  0.6535925189176744D-24
 0.630D+01   0.100D-05   0.1000000000000000D+01  0.6535925189007709D-23  0.1000000000000000D+01  0.6535925189007729D-23  0.1000000000000000D+01  0.6535925189007693D-23  0.1000000000000000D+01  0.6535925189007724D-23
 0.430D+00   0.100D-05   0.4568866945499436D+00  0.9378945443478894D-06  0.4568866945499436D+00  0.9378945443478894D-06  0.4568866945499436D+00  0.9378945443478894D-06  0.4568866945499435D+00  0.9378945443478893D-06
 0.130D+02   0.100D-01   0.1000000000000000D+01  0.4485365503481142D-75  0.1000000000000000D+01  0.4485365503481123D-75  0.1000000000000000D+01  0.4485365503481142D-75  0.1000000000000000D+01  0.4485365503481156D-75
 0.430D+01   0.100D+02 -0.9055021941027503D+34 -0.9379608607850867D+34 -0.9055021941027523D+34 -0.9379608607850810D+34 -0.9055021941027524D+34 -0.9379608607850814D+34 -0.9055021941027573D+34 -0.9379608607850862D+34
 0.630D-01   0.120D+02   0.1620124184557827D+62  0.1039653977037593D+61  0.1620124184557789D+62  0.1039653977037569D+61  0.1620124184557827D+62  0.1039653977037593D+61  0.1620124184557807D+62  0.1039653977037580D+61
 0.630D-01   0.150D+02   0.1857482893618806D+97 -0.6052285795414974D+96  0.1857482893618797D+97 -0.6052285795414943D+96  0.1857482893618806D+97 -0.6052285795414974D+96  0.1857482893618783D+97 -0.6052285795414903D+96
 0.630D-01   0.200D+02   0.8591804934512793+172 -0.1191477261295850+173  0.8591804934512993+172 -0.1191477261295879+173  0.8591804934512793+172 -0.1191477261295850+173  0.8591804934512906+172 -0.1191477261295866+173
 0.130D+02   0.160D+02  0.9749075032820455D+36  0.1347556133349446D+37  0.9749075032820759D+36  0.1347556133349415D+37  0.9749075032820455D+36  0.1347556133349446D+37  0.9749075032820488D+36  0.1347556133349445D+37
 0.100D+01   0.100D-01   0.8427423043912980D+00  0.4150936598121530D-02  0.8427423043912982D+00  0.4150936598121530D-02  0.8427423043912980D+00  0.4150936598121530D-02  0.8427423043912982D+00  0.4150936598121518D-02
 0.550D+00   0.100D-03   0.9999999999999926D+00  0.8223314414305329D-17  0.9999999999999927D+00  0.8223314414305349D-17  0.9999999999999926D+00  0.8223314414305329D-17  0.9999999999999927D+00  0.8223314414305334D-17
 0.190D+02   0.100D-01   0.1000000000000000D+01  0.4619188896571069-158  0.1000000000000000D+01  0.4619188896570986-158  0.1000000000000000D+01  0.4619188896571069-158  0.1000000000000000D+01  0.4619188896571069-158
 0.100D+01   0.280D+02  -Infinity                Infinity               -Infinity                Infinity               -Infinity                Infinity               -Infinity                NaN
 0.150D+02   0.150D+02   0.1000910969119025D+01 -0.2658046409880405D-01  0.1000910969119024D+01 -0.2658046409880407D-01  0.1000910969119025D+01 -0.2658046409880405D-01  0.1000910969119025D+01 -0.2658046409880409D-01
 0.260D+02   0.000D+00   0.1000000000000000D+01  0.0000000000000000D+00  0.1000000000000000D+01  0.0000000000000000D+00  0.1000000000000000D+01  0.0000000000000000D+00  0.1000000000000000D+01  0.0000000000000000D+00
 0.000D+00   0.266D+02   0.0000000000000000D+00  0.4132896053051959+306  0.0000000000000000D+00  0.4132896053051863+306  0.0000000000000000D+00  0.4132896053051739+306  0.0000000000000000D+00  0.4132896053051984+306
```

```
******************************************
******************************************
Table 8.  Comparison between values of erfi(z) calculated
from the present package and those calculated from
 Matlab (9.2.0.538062(R2017a)), Mathematica9 and Maple13
******************************************
    x         y         Re(Matlab)              Im(Matlab)             Re(Mathematica)         Im(Mathematica)         Re(Maple)               Im(Maple)              Re(present)             Im(present)
 0.630D-09  0.100D-09  0.7108788752701729D-09  0.1128379167095513D-09  0.7108788752701730D-09  0.1128379167095513D-09  0.7108788752701729D-09  0.1128379167095513D-09  0.7108788752701729D-09  0.1128379167095513D-09
 0.630D-08  0.100D-09  0.7108788752701730D-08  0.1128379167095513D-09  0.7108788752701730D-08  0.1128379167095513D-09  0.7108788752701729D-08  0.1128379167095513D-09  0.7108788752701730D-08  0.1128379167095513D-09
 0.230D-06  0.100D-09  0.2595272084319725D-06  0.1128379167095572D-09  0.2595272084319724D-06  0.1128379167095572D-09  0.2595272084319725D-06  0.1128379167095572D-09  0.2595272084319724D-06  0.1128379167095572D-09
 0.630D-01  0.100D-09  0.7118204889259457D-01  0.1132866603436267D-09  0.7118204889259455D-01  0.1132866603436267D-09  0.7118204889259457D-01  0.1132866603436267D-09  0.7118204889259457D-01  0.1132866603436267D-09
 0.630D+00  0.100D-08  0.8172721380624626D+00  0.1678133623772540D-08  0.8172721380624627D+00  0.1678133623772540D-08  0.8172721380624626D+00  0.1678133623772540D-08  0.8172721380624617D+00  0.1678133623772503D-08
 0.630D+01  0.100D-07  0.1566345974039032D+17  0.1948063216579133D+10  0.1566345974039038D+17  0.1948063217124518D+10  0.1566345974039036D+17  0.1948063216579138D+10  0.1566345974039032D+17  0.1948063216579134D+10
 0.230D+02  0.100D-04  0.1354844864027446D+229 0.6226384994616626+225  0.1354844864027350D+229 0.6226384994615598+225  0.1354844864027446D+229 0.6226384994616625+225  0.1354844864027385+229  0.6226384994616346+225
 0.510D-09  0.250D+00  0.5406072059818183D-09  0.2763263901682369D+00  0.5406072059818185D-09  0.2763263901682371D+00  0.5406072059818183D-09  0.2763263901682369D+00  0.5406072059818183D-09  0.2763263901682369D+00
 0.590D-09  0.150D+02  0.1279524608030563D-106 0.1000000000000000D+01  0.1279524608030580D-106 0.1000000000000000D+01  0.1279524608030563D-106 0.1000000000000000D+01  0.1279524608030563D-106 0.1000000000000000D+01
 0.630D-09  0.250D+02  0.2616642670093135D-280 0.1000000000000000D+01  0.1856632045757707D-302 0.1000000000000000D+01  0.2616642670093135D-280 0.1000000000000000D+01  0.2616642670093134D-280 0.1000000000000000D+01
 0.630D-07  0.100D-09  0.7108788752701738D-07  0.1128379167095517D-09  0.7108788752701738D-07  0.1128379167095517D-09  0.7108788752701740D-07  0.1128379167095517D-09  0.7108788752701738D-07  0.1128379167095517D-09
 0.630D-07  0.100D-06  0.7108788752701667D-07  0.1128379167095513D-06  0.7108788752701667D-07  0.1128379167095513D-06  0.7108788752701668D-07  0.1128379167095513D-06  0.7108788752701667D-07  0.1128379167095513D-06
 0.630D+01  0.100D-07  0.1566345974039044D+17  0.1948063216579138D+08  0.1566345974039051D+17  0.1948062894718035D+08  0.1566345974039048D+17  0.1948063216579143D+08  0.1566345974039044D+17  0.1948063216579139D+08
 0.630D+01  0.100D-07  0.1566345974039032D+17  0.1948063216579133D+10  0.1566345974039038D+17  0.1948063217124518D+10  0.1566345974039036D+17  0.1948063216579138D+10  0.1566345974039032D+17  0.1948063216579134D+10
 0.630D+01  0.100D-03  0.1566344746759387D+17  0.1948062694628111D+14  0.1566344746759394D+17  0.1948062694628334D+14  0.1566344746759390D+17  0.1948062694628116D+14  0.1566344746759385D+17  0.1948062694628111D+14
 0.630D+01  0.100D+01  0.5638834068500719D+16 -0.7247961961672762D+15  0.5638834068500730D+16 -0.7247961961672785D+15  0.5638834068500733D+16 -0.7247961961672758D+15  0.5638834068500740D+16 -0.7247961961672788D+15
 0.630D+01  0.250D+02  0.1230377519082839D-255 0.1000000000000000D+01  0.1230377519082947D-255 0.1000000000000000D+01  0.1230377519082849D-255 0.1000000000000000D+01  0.1230377519082916D-255 0.1000000000000000D+01
 0.630D+01  0.260D+02  0.8519503654259445D-278 0.1000000000000000D+01  0.8519503654259705D-278 0.1000000000000000D+01  0.8519503654259511D-278 0.1000000000000000D+01  0.8519503654259835D-278 0.1000000000000000D+01
 0.630D+01  0.120D+02  0.1607180360354913D-46  0.1000000000000000D+01  0.1607180360354879D-46  0.1000000000000000D+01  0.1607180360354923D-46  0.1000000000000000D+01  0.1607180360354905D-46  0.1000000000000000D+01
 0.130D+02  0.150D+02  0.1243071659465133D-25  0.1000000000000000D+01  0.1243071659465130D-25  0.1000000000000000D+01  0.1243071659465133D-25  0.1000000000000000D+01  0.1243071659465133D-25  0.1000000000000000D+01
 0.130D+01  0.200D+02  0.2842335234224522D-174 0.1000000000000000D+01  0.2842335234224539D-174 0.1000000000000000D+01  0.2842335234224522D-174 0.1000000000000000D+01  0.2842335234224511D-174 0.1000000000000000D+01
 0.230D+01  0.100D+02  0.3164123215495136D-42  0.1000000000000000D+01  0.3164123215495133D-42  0.1000000000000000D+01  0.3164123215495129D-42  0.1000000000000000D+01  0.3164123215495102D-42  0.1000000000000000D+01
 0.260D+02  0.100D-01  0.7216470199790084D+292 0.4128185256330935+292  0.7216470199790552D+292 0.4128185256331204+292  0.7216470199790084D+292 0.4128185256330935+292  0.7216470199790265D+292 0.4128185256331038+292
 0.730D+01  0.100D-03  0.1086000033675769D+23  0.1570391168396834D+20  0.1086000033675770D+23  0.1570391168396819D+20  0.1086000033675772D+23  0.1570391168396838D+20  0.1086000033675770D+23  0.1570391168396836D+20
 0.630D+01  0.100D-06  0.1566345974037818D+17  0.1948063216578616D+11  0.1566345974037828D+17  0.1948063216608607D+11  0.1566345974037821D+17  0.1948063216578621D+11  0.1566345974037822D+17  0.1948063216578622D+11
 0.630D+01  0.100D-05  0.1566345973916317D+17  0.1948063216526943D+12  0.1566345973916319D+17  0.1948063216505383D+12  0.1566345973916320D+17  0.1948063216526948D+12  0.1566345973916313D+17  0.1948063216526940D+12
 0.430D+00  0.100D-05  0.5168422696512865D+00  0.1357550859430707D-05  0.5168422696512868D+00  0.1357550859430707D-05  0.5168422696512865D+00  0.1357550859430707D-05  0.5168422696512867D+00  0.1357550859430707D-05
 0.130D+01  0.100D-01  0.1046506224026099D+73  0.2775262271063284D+72  0.1046506224026087D+73  0.2775262271063254D+72  0.1046506224026099D+73  0.2775262271063284D+72  0.1046506224026095D+73  0.2775262271063274D+72
 0.430D+01  0.100D+02 -0.2060444288437954D-36  0.1000000000000000D+01 -0.2060444288437958D-36  0.1000000000000000D+01 -0.2060444288437954D-36  0.1000000000000000D+01 -0.2060444288437947D-36  0.1000000000000000D+01
 0.630D-01  0.120D+02  0.1359682283481483D-63  0.1000000000000000D+01  0.1359682283481485D-63  0.1000000000000000D+01  0.1359682283481483D-63  0.1000000000000000D+01  0.1359682283481500D-63  0.1000000000000000D+01
 0.630D-01  0.150D+02  0.6866248044893615D-99  0.1000000000000000D+01  0.6866248044893749D-99  0.1000000000000000D+01  0.6866248044893615D-99  0.1000000000000000D+01  0.6866248044893696D-99  0.1000000000000000D+01
 0.630D-01  0.200D+02  0.3140804829211560-175  0.1000000000000000D+01  0.3140804829211534-175  0.1000000000000000D+01  0.3140804829211560-175  0.1000000000000000D+01  0.3140804829211500-175  0.1000000000000000D+01
 0.130D+02  0.160D+02  0.4111408922217984D-39  0.1000000000000000D+01  0.4111408922217920D-39  0.1000000000000000D+01  0.4111408922217984D-39  0.1000000000000000D+01  0.4111408922217990D-39  0.1000000000000000D+01
 0.100D+01  0.100D-01  0.1650119059098099D+01  0.3066945879694022D-01  0.1650119059098099D+01  0.3066945879694022D-01  0.1650119059098099D+01  0.3066945879694022D-01  0.1650119059098099D+01  0.3066945879694027D-01
 0.550D+01  0.100D-03  0.1432098320459191D+13  0.1548328285952877D+10  0.1432098320459185D+13  0.1548328285952778D+10  0.1432098320459191D+13  0.1548328285952877D+10  0.1432098320459190D+13  0.1548328285952876D+10
 0.190D+02  0.100D+01  0.6385144550807778+155  0.1617880125819369+155  0.6385144550807715+155  0.1617880125819336+155  0.6385144550807778+155  0.1617880125819369+155  0.6385144550807778+155  0.1617880125819369+155
 0.100D+01  0.280D+02  0.0000000000000000D+00  0.1000000000000000D+01  0.0000000000000000D+00  0.1000000000000000D+01  0.0000000000000000D+00  0.1000000000000000D+01  0.0000000000000000D+00  0.1000000000000000D+01
 0.150D+02  0.150D+02 -0.2658046409880405D-01  0.1000910969119025D+01 -0.2658046409880407D-01  0.1000910969119024D+01 -0.2658046409880405D-01  0.1000910969119025D+01 -0.2658046409880409D-01  0.1000910969119025D+01
 0.260D+02  0.000D+00  0.8314637164730988D+292 0.0000000000000000D+00  0.8314637164730799D+292 0.0000000000000000D+00  0.8314637164730988D+292 0.0000000000000000D+00  0.8314637164730988D+292 0.0000000000000000D+00
 0.000D+00  0.266D+02  0.0000000000000000D+00  0.1000000000000000D+01  0.0000000000000000D+00  0.1000000000000000D+01  0.0000000000000000D+00  0.1000000000000000D+01  0.0000000000000000D+00  0.1000000000000000D+01
```

```
******************************************
******************************************
Table 9.   Comparison between values of Daw(z) calculated
from the present package and those calculated from
 Matlab (9.2.0.538062(R2017a)), Mathematica9 and Maple3
******************************************
    x            y          Re(Matlab)                Im(Matlab)              Re(Mathematica)           Im(Mathematica)            Re(Maple)                  Im(Maple)                 Re(present)                Im(present)
 0.630D-09  0.100D-09   0.630000000000000D-09   0.100000000000000D-09   0.630000000000000D-09   0.100000000000000D-09   0.630000000000000D-09   0.100000000000000D-09   0.630000000000000D-09   0.100000000000000D-09
 0.630D-08  0.100D-09   0.630000000000000D-08   0.999999999999999D-10   0.630000000000000D-08   0.999999999999999D-10   0.630000000000000D-08   0.999999999999999D-10   0.630000000000000D-08   0.999999999999999D-10
 0.230D-06  0.100D-09   0.229999999999919D-06   0.999999999998942D-10   0.229999999999919D-06   0.999999999998942D-10   0.229999999999919D-06   0.999999999998942D-10   0.229999999999919D-06   0.999999999998942D-10
 0.630D-01  0.100D-09   0.6283356634989649D-01  0.9920829706399130D-10  0.6283356634989649D-01  0.9920829706399130D-10  0.6283356634989649D-01  0.9920829706399130D-10  0.6283356634989649D-01  0.9920829706399130D-10
 0.630D+00  0.100D-08   0.4870125516138508D+00  0.3863641849665480D-09  0.4870125516138508D+00  0.3863641849665479D-09  0.4870125516138508D+00  0.3863641849665479D-09  0.4870125516138508D+00  0.3863641849665480D-09
 0.630D+01  0.100D-07   0.8040529489538834D-01 -0.1310671568189310D-09  0.8040529489538865D-01 -0.1310671540193449D-09  0.8040529489538834D-01 -0.1310671568189310D-09  0.8040529489538835D-01 -0.1310671568189315D-09
 0.230D+02  0.100D-04   0.2175973635712256D-01 -0.9478724278279226D-08  0.2175973635712200D-01 -0.9478724279222841D-08  0.2175973635712256D-01 -0.9478724278279226D-08  0.2175973635712256D-01 -0.9478724278279180D-08
 0.510D-09  0.250D+00   0.5764738587346685D-09  0.2606817989594842D+00  0.5764738587346688D-09  0.2606817989594844D+00  0.5764738587346685D-09  0.2606817989594842D+00  0.5764738587346686D-09  0.2606817989594842D+00
 0.590D-09  0.150D+02   0.8161624977321701D+90  0.4611087557808870D+98  0.8161624977321701D+90  0.4611087557808870D+98  0.8161624977321701D+90  0.4611087557808870D+98  0.8161624977321701D+90  0.4611087557808870D+98
 0.630D-09  0.250D+02   0.7584145984783490D+264 0.2407665391994758D+272 0.1111625755203242D+287 0.3393241010998907D+294 0.7584145984783490D+264 0.2407665391994758D+272 0.7584145984783488D+264 0.2407665391994758D+272
 0.630D-07  0.100D-09   0.629999999999984D-07   0.999999999999922D-10   0.629999999999982D-07   0.999999999999920D-10   0.629999999999984D-07   0.999999999999922D-10   0.629999999999982D-07   0.999999999999922D-10
 0.630D-07  0.100D-06   0.630000000000109D-07   0.999999999999986D-07   0.630000000000107D-07   0.999999999999986D-07   0.630000000000109D-07   0.999999999999988D-07   0.630000000000108D-07   0.999999999999986D-07
 0.630D+01  0.100D-09   0.8040529489538834D-01 -0.1310671568189310D-11  0.8040529489538865D-01 -0.1310688090297495D-11  0.8040529489538834D-01 -0.1310671568189310D-11  0.8040529489538835D-01 -0.1310671568189315D-11
 0.630D+01  0.100D-07   0.8040529489538834D-01 -0.1310671568189310D-09  0.8040529489538865D-01 -0.1310671540193449D-09  0.8040529489538834D-01 -0.1310671568189310D-09  0.8040529489538835D-01 -0.1310671568189315D-09
 0.630D+01  0.100D-03   0.8040529487371820D-01 -0.1310671567825662D-05  0.8040529487371867D-01 -0.1310671567814718D-05  0.8040529487371820D-01 -0.1310671567825662D-05  0.8040529487371822D-01 -0.1310671567825666D-05
 0.630D+01  0.100D+01   0.7829843135459250D-01 -0.1275347937863088D-01  0.7829843135459262D-01 -0.1275347937863092D-01  0.7829843135459250D-01 -0.1275347937863088D-01  0.7829843135459279D-01 -0.1275347937863093D-01
 0.630D+01  0.250D+02   0.1039159219998316D+255 0.9300772466823997D+254 0.1039159219998266D+255 0.9300772466823360D+254 0.1039159219998322D+255 0.9300772466823869D+254 0.1039159219998265D+255 0.9300772466823363D+254
 0.630D+01  0.260D+02   0.1507783230583405D+277 0.1260820133420931D+277 0.1507783230583289D+277 0.1260820133420895D+277 0.1507783230583414D+277 0.1260820133420912D+277 0.1507783230583289D+277 0.1260820133420895D+277
 0.630D+01  0.120D+02   0.6963860290336417D+45  0.1630929480225655D+46  0.6963860290336307D+45  0.1630929480225659D+46  0.6963860290336477D+45  0.1630929480225648D+46  0.6963860290336308D+45  0.1630929480225659D+46
 0.130D+02  0.150D+02   0.7937662726696219D+24  0.1675136726074209D+25  0.7937662726696219D+24  0.1675136726074209D+25  0.7937662726696219D+24  0.1675136726074209D+25  0.7937662726696219D+24  0.1675136726074209D+25
 0.130D+02  0.200D+02   0.8424289527571354D+173 -0.1391691798227752D+173 0.8424289527571373D+173 -0.1391691798227755D+173 0.8424289527571354D+173 -0.1391691798227752D+173 0.8424289527571374D+173 -0.1391691798227755D+173
 0.230D+01  0.100D+02   0.1083128394332845D+42  -0.5190843336239377D+41  0.1083128394332851D+42  -0.5190843336239456D+41  0.1083128394332842D+42  -0.5190843336239416D+41  0.1083128394332851D+42  -0.5190843336239457D+41
 0.260D+02  0.100D-01   0.1924502199435739D-01  -0.7412921854824288D-05  0.1924502199435815D-01  -0.7412921854822422D-05  0.1924502199435739D-01  -0.7412921854824289D-05  0.1924502199435739D-01  -0.7412921854824279D-05
 0.730D+01  0.100D-03   0.6915479482198875D-01  -0.9660004598137924D-06  0.6915479482198875D-01  -0.9660004598124207D-06  0.6915479482198875D-01  -0.9660004598124207D-06  0.6915479482198875D-01  -0.9660004598124298D-06
 0.630D+01  0.100D-06   0.8040529489538832D-01  -0.1310671568189310D-08  0.8040529489538863D-01  -0.1310671566650571D-08  0.8040529489538832D-01  -0.1310671568189310D-08  0.8040529489538832D-01  -0.1310671568189314D-08
 0.630D+01  0.100D-05   0.8040529489538618D-01  -0.1310671568189274D-07  0.8040529489538646D-01  -0.1310671569296230D-07  0.8040529489538617D-01  -0.1310671568189273D-07  0.8040529489538618D-01  -0.1310671568189314D-07
 0.430D+00  0.100D-05   0.3807166899583291D+00   0.6725836466366695D-06  0.3807166899583292D+00   0.6725836466366695D-06  0.3807166899583291D+00   0.6725836466366695D-06  0.3807166899583291D+00   0.6725836466366693D-06
 0.130D+02  0.100D-01   0.3857633206650098D-01  -0.2985234384101036D-04  0.3857633206650069D-01  -0.2985234384100776D-04  0.3857633206650098D-01  -0.2985234384101037D-04  0.3857633206650097D-01  -0.2985234384101034D-04
 0.430D+01  0.100D+02  -0.2052605644382073D+36  -0.8528608909254047D+35 -0.2052605644382083D+36  -0.8528608909253993D+35 -0.2052605644382073D+36  -0.8528608909253949D+35 -0.2052605644382083D+36  -0.8528608909253993D+35
 0.630D+01  0.120D+02   0.3044216301796831D+63   0.1791952772920917D+62  0.3044216301796794D+63   0.1791952772920895D+62  0.3044216301796831D+63   0.1791952772920917D+62  0.3044216301796794D+63   0.1791952772920895D+62
 0.630D+01  0.150D+02   0.4360818821725557D+98  -0.1441276173418956D+98  0.4360818821725504D+98  -0.1441276173418939D+98  0.4360818821725557D+98  -0.1441276173418956D+98  0.4360818821725505D+98  -0.1441276173418939D+98
 0.630D+01  0.200D+02   0.2684006948997963D+174 -0.3746958739955385D+174 0.2684006948998007D+174 -0.3746958739955445D+174 0.2684006948997963D+174 -0.3746958739955385D+174 0.2684006948998007D+174 -0.3746958739955446D+174
 0.130D+02  0.160D+02   0.5202336425049914D+38   0.1389653791837263D+38  0.5202336425049913D+38   0.1389653791837262D+38  0.5202336425049914D+38   0.1389653791837263D+38  0.5202336425049914D+38   0.1389653791837262D+38
 0.100D+01  0.100D-01   0.5381256994780437D+00  -0.7619488743398086D-03  0.5381256994780438D+00  -0.7619488743398104D-03  0.5381256994780437D+00  -0.7619488743398086D-03  0.5381256994780438D+00  -0.7619488743397998D-03
 0.550D+01  0.100D-05   0.9249323227728154D-01  -0.1742555541173757D-05  0.9249323227728121D-01  -0.1742555541179759D-05  0.9249323227728154D-01  -0.1742555541173757D-05  0.9249323227728153D-01  -0.1742555541173745D-05
 0.190D+02  0.100D+01   0.2627908698045743D-01  -0.1386957125910295D-02  0.2627908698045715D-01  -0.1386957125910346D-02  0.2627908698045743D-01  -0.1386957125910295D-02  0.2627908698045743D-01  -0.1386957125910295D-02
 0.100D+01  0.280D+02  -Infinity                 Infinity               -Infinity                 Infinity               -Infinity                 Infinity               -Infinity                 Infinity
 0.150D+02  0.150D+02  -0.5888963480848108D+00  -0.6637663396724557D+00 -0.5888963480848102D+00  -0.6637663396724550D+00 -0.5888963480848108D+00  -0.6637663396724557D+00 -0.5888963480848108D+00  -0.6637663396724557D+00
 0.260D+02  0.000D+00   0.1924502485184063D-01   0.0000000000000000D+00  0.1924502485184020D-01   0.0000000000000000D+00  0.1924502485184063D-01   0.0000000000000000D+00  0.1924502485184064D-01   0.0000000000000000D+00
 0.000D+00  0.266D+02   0.0000000000000000D+00   0.1725633471960343D+308 0.0000000000000000D+00   0.1725633471960353D+308 0.0000000000000000D+00   0.1725633471960251D+308 0.0000000000000000D+00   0.1725633471960353D+308
```

```
******************************************
******************************************
Table 10. Comparison between values of FresnelS(z) calculated
from the present package and those calculated from
 Matlab (9.2.0.538062(R2017a)), Mathematica9 and Maple13
******************************************
    x           y           Re(Matlab)              Im(Matlab)              Re(Mathematica)         Im(Mathematica)         Re(Maple)               Im(Maple)               Re(present)             Im(present)
 0.630D-09  0.100D-09   0.1210282861832200D-27   0.6182130743489115D-28   0.1276256307557585D-27   0.2025803662789818D-28   0.1210282861832200D-27   0.6182130743489116D-28   0.1210282861832200D-27   0.6182130743489115D-28
 0.630D-08  0.100D-09   0.1308253428734398D-24   0.6233967022273347D-26   0.1308913163191651D-24   0.2077639941574050D-26   0.1308253428734398D-24   0.6233967022273348D-26   0.1308253428734397D-24   0.6233967022273346D-26
 0.230D-06  0.100D-09   0.6370622689872950D-20   0.8309512045146226D-23   0.6370625098427318D-20   0.2769836999316225D-23   0.6370622689872952D-20   0.8309512045146229D-23   0.6370622689872949D-20   0.8309512045146226D-23
 0.630D-01  0.100D-09   0.1309239395510383D-03   0.6234450233189308D-12   0.1309239395510383D-03   0.6234513030616896D-12   0.1309239395510383D-03   0.6234450233189760D-12   0.1309239395510383D-03   0.6234450233189758D-12
 0.630D+00  0.100D-08   0.1273340391859734D+00   0.5838388163123306D-09   0.1273340391859734D+00   0.5838388570588143D-09   0.1273340391859735D+00   0.5838388163123305D-09   0.1273340391859734D+00   0.5838388163123302D-09
 0.630D+01  0.100D-07   0.4555454305043985D+00  -0.4679298142606107D-10   0.4555454305043977D+00  -0.4679294306591841D-10   0.4555454305043987D+00  -0.4679298142605792D-10   0.4555454305043853D+00  -0.4679298142605787D-10
 0.230D+02  0.100D-04   0.4999916725048591D+00   0.1000000087017015D-04   0.4999916725048578D+00   0.1000000087014726D-04   0.4999916725048605D+00   0.1000000087017014D-04   0.4999916725048115D+00   0.1000000087017012D-04
 0.510D-09  0.250D+00  -0.4998874156807607D-10  -0.8175600235777757D-02  -0.4998871773360066D-10  -0.8175600235777755D-02  -0.4998874156807591D-10  -0.8175600235777758D-02  -0.4998874156807590D-10  -0.8175600235777755D-02
 0.590D-09  0.150D+02  -0.5900000000000077D-09  -0.4996997980970270D+00  -0.5900000987429972D-09  -0.4996997980970240D+00  -0.5900000000000000D-09  -0.4996997980970300D+00  -0.5899999999999990D-09  -0.4996997980972190D+00
 0.630D-09  0.250D+00  -0.6299999999999798D-09  -0.4999935154694762D+00  -0.3973874519835233D-16  -0.4877573202131747D-09  -0.6300000000000280D-09  -0.4999935154694769D+00  -0.6299999999999980D-09  -0.4999935154692640D+00
 0.630D-07  0.100D-09   0.1309233134403419D-21   0.6234485385061188D-24   0.1309239731747992D-21   0.2078158304361892D-24   0.1309233134403420D-21   0.6234485385061190D-24   0.1309233134403419D-21   0.6234485385061188D-24
 0.630D-07  0.100D-06  -0.8586773828387569D-21   0.9985028650659554D-22  -0.8595256252368686D-21   0.1003846792174059D-21  -0.8586773828387573D-21   0.9985028650659562D-22  -0.8586773828387569D-21   0.9985028650659559D-22
 0.630D+01  0.100D-09   0.4555454305043985D+00  -0.4679298142606107D-10   0.4555454305043977D+00  -0.4679294306591841D-10   0.4555454305043987D+00  -0.4679298142605792D-10   0.4555454305043853D+00  -0.4679298142605787D-10
 0.630D+01  0.100D-07   0.4555454305043976D+00  -0.4679298142605799D-08   0.4555454305043977D+00  -0.4679298109799592D-08   0.4555454305043978D+00  -0.4679298142605818D-08   0.4555454305043853D+00  -0.4679298142605817D-08
 0.630D+01  0.100D-03   0.4555453430467746D+00  -0.4679301243873927D-04   0.4555453430467747D+00  -0.4679301243876245D-04   0.4555453430467746D+00  -0.4679301243873937D-04   0.4555453430467744D+00  -0.4679301243865064D-04
 0.630D+01  0.100D+01   0.3259038775999915D+07  -0.9300208548761779D+07   0.3259038775999944D+07  -0.9300208548761757D+07   0.3259038775999923D+07  -0.9300208548761779D+07   0.3259038775999811D+07  -0.9300208548761800D+07
 0.630D+01  0.250D+02  -0.3554860935704188D+213 -0.3203042440857180D+213 -0.3554860935704541D+213 -0.3203042440856957D+213 -0.3554860935703959D+213 -0.3203042440857403D+213 -0.3554860935704820D+213 -0.3203042440858674D+213
 0.630D+01  0.260D+02  -0.1204038050090805D+222  0.1361012791371615D+222 -0.1204038050090677D+222  0.1361012791371676D+222 -0.1204038050090854D+222  0.1361012791371530D+222 -0.1204038050090850D+222  0.1361012791371248D+222
 0.630D+01  0.120D+02  -0.1359772497374774D+102  0.9330549183612566D+101 -0.1359772497374750D+102  0.9330549183612255D+101 -0.1359772497374827D+102  0.9330549183612566D+101 -0.1359772497374866D+102  0.9330549183613046D+101
 0.130D+02  0.150D+02  -0.5942189069025935D+264  0.6857969362553382D+264 -0.5942189069025582D+264  0.6857969362553023D+264 -0.5942189069025940D+264  0.6857969362553339D+264 -0.5942189069026832D+264  0.6857969362554658D+264
 0.130D+02  0.200D+02   0.1237947427193671D+34  -0.2014306270028455D+34   0.1237947427193599D+34  -0.2014306270028513D+34   0.1237947427193716D+34  -0.2014306270028460D+34   0.1237947427193897D+34  -0.2014306270028485D+34
 0.230D+01  0.100D+02   0.3631076657574397D+30  -0.8586224684176883D+29   0.3631076657574405D+30  -0.8586224684176507D+29   0.3631076657574440D+30  -0.8586224684177132D+29   0.3631076657574456D+30  -0.8586224684176628D+29
 0.260D+02  0.100D-01   0.4834410699595388D+00  -0.6326347012213576D-06   0.4834410699595387D+00  -0.6326347001719534D-06   0.4834410699595388D+00  -0.6326347006416399D-06   0.4834410699595386D+00  -0.6326347001262000D-06
 0.730D+01  0.100D-03   0.5189473783051683D+00   0.8980283652621740D-04   0.5189473783051683D+00   0.8980283652625001D-04   0.5189473783051680D+00   0.8980283652621786D-04   0.5189473783051683D+00   0.8980283652659793D-04
 0.630D+01  0.100D-06   0.4555454305043110D+00  -0.4679298142608870D-07   0.4555454305043111D+00  -0.4679298142750D-07      0.4555454305043108D+00  -0.4679298142608786D-07   0.4555454305043470D+00  -0.4679298142608888D-07
 0.630D+01  0.100D-05   0.4555454304956527D+00  -0.4679298142915895D-06   0.4555454304956528D+00  -0.4679298143474182D-06   0.4555454304956529D+00  -0.4679298142915893D-06   0.4555454305131301D+00  -0.4679298142915913D-06
 0.430D+00  0.100D-01   0.4137960430796478D-01   0.2863740555405401992D-06  0.4137960430796477D-01   0.2863740555059D-06      0.4137960430796480D-01   0.2863740555405401999D-06  0.4137960430796477D-01   0.2863740554504199D-06
 0.130D+02  0.100D-01   0.4999537211098605D-01   0.1028032147726659D-01   0.4999537211098603D-01   0.1028032147726653D-01   0.4999537211098607D-01   0.1028032147726659D-01   0.4999537211098617D-01   0.1028032147726610D-01
 0.430D+01  0.100D+02  -0.2438578235593966D+57  -0.6372912157595411D+57  -0.2438578235593918D+57  -0.6372912157595408D+57  -0.2438578235594011D+57  -0.6372912157595420D+57  -0.2438578235593919D+57  -0.6372912157595411D+57
 0.630D-01  0.120D+02   0.4386953709430511D-03  -0.3561716912083346D+00   0.4386953709470129D-03  -0.3561716912083345D+00   0.4386953709420732D-03  -0.3561716912083348D+00   0.4386953709430291D-03  -0.3561716912083304D+00
 0.630D-01  0.150D+02  -0.2060221319274167D+00  -0.4992810486477227D+00  -0.2060221319274166D+00  -0.4992810486477147D+00  -0.2060221319274166D+00  -0.4992810486477215D+00  -0.2060221319274229D+00  -0.4992810486477265D+00
 0.630D-01  0.200D+02   0.1615654480771972D-02  -0.8307427050562236D-01   0.1615654480778416D-02  -0.8307427050562262D-01   0.1615654480766620D-02  -0.8307427050562211D-01   0.1615654480731143D-02  -0.8307427050559757D-01
 0.130D+02  0.160D+02   0.3700036576630200D+282  0.3005334592424596D+282  0.3700036576630358D+282  0.3005334592424637D+282  0.3700036576630594D+282  0.3005334592424825D+282  0.3700036576630693D+282  0.3005334592425084D+282
 0.100D+01  0.100D-01   0.4382591350519963D+00   0.1000164499056029D-01   0.4382591350519964D+00   0.1000164499056032D-01   0.4382591350519963D+00   0.1000164499056029D-01   0.4382591350519965D+00   0.1000164499056022D-01
 0.550D+01  0.100D-03   0.5536841425965034D+00  -0.3826836179481189D-04   0.5536841425965037D+00  -0.3826836179478648D-04   0.5536841425965039D+00  -0.3826836179481149D-04   0.5536841425965033D+00  -0.3826836179476878D-04
 0.190D+02  0.100D+01  -0.6999628325952581D+24   0.3622612080927746D+23  -0.6999628325952543D+24   0.3622612080927184D+23  -0.6999628325952504D+24   0.3622612080927570D+23  -0.6999628325952505D+24   0.3622612080935166D+23
 0.100D+01  0.280D+02   0.9050342989617665D+36   0.3195609198840348D+35   0.9050342989617654D+36   0.3195609198833903D+35   0.9050342989617598D+36   0.3195609198840666D+35   0.9050342989616937D+36   0.3195609198854360D+35
 0.150D+02  0.150D+02  -0.5124909928846749D+305  0.5124909928846749D+305 -0.5124909928846555D+305  0.5124909928846555D+305 -0.5124909928846823D+305  0.5124909928846823D+305 -0.5124909928846552D+305  0.5124909928846552D+305
 0.260D+02  0.000D+00   0.4877573202131747D+00   0.0000000000000000D+00   0.4877573202131748D+00   0.0000000000000000D+00   0.4877573202131747D+00   0.0000000000000000D+00   0.4877573202129533D+00   0.0000000000000000D+00
 0.000D+00  0.266D+02   0.0000000000000000D+00  -0.4907830617995415D+00   0.0000000000000000D+00  -0.4907830617995422D+00   0.0000000000000000D+00  -0.4907830617995416D+00   0.0000000000000000D+00  -0.4907830617993149D+00
```

```
******************************************
******************************************
Table 11. Comparison between values of FresnelC(z) calculated
from the present package and those calculated from
 Matlab (9.2.0.538062(R2017a)), Mathematica9 and Maple13
******************************************
     x         y         Re(Matlab)              Im(Matlab)             Re(Mathematica)         Im(Mathematica)         Re(Maple)               Im(Maple)               Re(present)             Im(present)
 0.630D-09  0.100D-09  0.630000000000000D-09   0.100000000000000D-09   0.630000000000000D-09   0.100000000000000D-09   0.630000000000000D-09   0.100000000000000D-09   0.630000000000000D-09   0.100000000000000D-09
 0.630D-08  0.100D-09  0.630000000000000D-08   0.100000000000000D-09   0.630000000000000D-08   0.100000000000000D-09   0.630000000000000D-08   0.100000000000000D-09   0.630000000000000D-08   0.100000000000000D-09
 0.230D-06  0.100D-09  0.230000000000000D-06   0.100000000000000D-09   0.230000000000000D-06   0.100000000000000D-09   0.230000000000000D-06   0.100000000000000D-09   0.230000000000000D-06   0.100000000000000D-09
 0.630D-01  0.100D-09  0.629975512653883D-01   0.999805656262978D-10   0.629975512653883D-01   0.999805655954585D-10   0.629975512653883D-01   0.999805656262975D-10   0.629975512653884D-01   0.999805656262975D-10
 0.630D+00  0.100D-08  0.605949325118742D+00   0.811869593325809D-09   0.605949325118742D+00   0.811869596950979D-09   0.605949325118742D+00   0.811869593325810D-09   0.605949325118743D+00   0.811869593325810D-09
 0.630D+01  0.100D-07  0.476004455353067D+00   0.883765630088697D-08   0.476004455353067D+00   0.883765633432947D-08   0.476004455353067D+00   0.883765630088696D-08   0.476004455353081D+00   0.883765630088696D-08
 0.230D+02  0.100D-04  0.513839548849339D+00   0.523598857539315D-15   0.513839548849339D+00   0.516182883908395D-15   0.513839548849339D+00   0.522685504000000D-15   0.513839548849387D+00   0.523598857609900D-15
 0.510D-09  0.250D+02  0.507544210602820D-09   0.249759150356543D+00   0.507544210602824D-09   0.249759150356543D+00   0.507544210602820D-09   0.249759150356543D+00   0.507544210602820D-09   0.249759150356543D+00
 0.590D-09  0.150D+02  0.827180612553027D-24   0.521220531674373D+00  -0.721546722921759D-16   0.521220531674373D+00  -0.415930000000000D-23   0.521220531674373D+00   0.107536192933603D-27   0.521220531674354D+00
 0.630D-09  0.250D+02  0.000000000000000D+00   0.512732385539770D+00   0.629999977893479D-09   0.499994235272720D+00  -0.392622000000000D-22   0.512732385539770D+00   0.130924303042027D-27   0.512732385539982D+00
 0.630D-07  0.100D-06  0.630000000000000D-07   0.100000000000000D-06   0.630000000000000D-07   0.100000000000000D-06   0.630000000000000D-07   0.100000000000000D-06   0.630000000000000D-07   0.100000000000000D-06
 0.630D+01  0.100D-07  0.476004455353067D+00   0.883765630088748D-10   0.476004455353067D+00   0.883765449829814D-10   0.476004455353067D+00   0.883765630088903D-10   0.476004455353082D+00   0.883765630088690D-10
 0.630D+01  0.100D-07  0.476004455353067D+00   0.883765630088697D-08   0.476004455353067D+00   0.883765633432947D-08   0.476004455353067D+00   0.883765630088696D-08   0.476004455353081D+00   0.883765630088696D-08
 0.630D+01  0.100D-03  0.476004409046638D+00   0.883766204626621D-04   0.476004409046639D+00   0.883766204625727D-04   0.476004409046638D+00   0.883766204625727D-04   0.476004409046652D+00   0.883766204626620D-04
 0.630D+01  0.100D+01 -0.930020804876177D+07  -0.325903827599991D+07  -0.930020804876175D+07  -0.325903827599994D+07  -0.930020804876177D+07  -0.325903827599992D+07  -0.930020804876180D+07  -0.325903827599981D+07
 0.630D+01  0.250D+02 -0.320304244085718D+213  0.355486093570418D+213 -0.320304244085695D+213  0.355486093570454D+213 -0.320304244085740D+213  0.355486093570395D+213 -0.320304244085867D+213  0.355486093570482D+213
 0.630D+01  0.260D+02  0.136101279137161D+222  0.120403805009080D+222  0.136101279137167D+222  0.120403805009067D+222  0.136101279137153D+222  0.120403805009085D+222  0.136101279137124D+222  0.120403805009085D+222
 0.630D+01  0.120D+02  0.933054918361256D+101  0.135977249737477D+102  0.933054918361225D+101  0.135977249737475D+102  0.933054918361256D+101  0.135977249737482D+102  0.933054918361304D+101  0.135977249737486D+102
 0.130D+02  0.150D+02  0.685796936255338D+264  0.594218906902558D+264  0.685796936255302D+264  0.594218906902558D+264  0.685796936255333D+264  0.594218906902594D+264  0.685796936255466D+264  0.594218906902683D+264
 0.130D+02  0.200D+02 -0.201430627002845D+34  -0.123794742719367D+34  -0.201430627002851D+34  -0.123794742719359D+34  -0.201430627002846D+34  -0.123794742719371D+34  -0.201430627002848D+34  -0.123794742719389D+34
 0.230D+01  0.100D+02 -0.858622468417688D+29  -0.363107665757439D+30  -0.858622468417650D+29  -0.363107665757440D+30  -0.858622468417713D+29  -0.363107665757443D+30  -0.858622468417663D+29  -0.363107665757445D+30
 0.260D+02  0.100D-01  0.499993890129386D+00   0.111496648351478D-01   0.499993890129390D+00   0.111496648351478D-01   0.499993890129389D+00   0.111496648351478D-01   0.499993890129395D+00   0.111496648351480D-01
 0.730D-03  0.100D-03  0.539268118633652D+00  -0.439939550798181D-04   0.539268118633652D+00  -0.439939550798443D-04   0.539268118633652D+00  -0.439939550798172D-04   0.539268118633652D+00  -0.439939550798185D-04
 0.630D+01  0.100D-06  0.476004455353021D-07   0.883765630372868D-07   0.476004455353021D-07   0.883765630372887D-07   0.476004455353022D-07   0.883765630372870D-07   0.476004455353035D-07   0.883765630372865D-07
 0.630D+01  0.100D-05  0.476004455484373D+00   0.883765630146145D-05   0.476004455484372D+00   0.883765630153605D-05   0.476004455484373D+00   0.883765630146146D-05   0.476004455484513D+00   0.883765630146144D-05
 0.430D+00  0.100D-05  0.426386850340233D+00   0.958117894814321D-05   0.426386850340233D+00   0.958117894813441D-05   0.426386850340233D+00   0.958117894814321D-05   0.426386850340232D+00   0.958117894814321D-05
 0.130D+01  0.100D-01  0.526555692463648D+00   0.550060331870065D-02   0.526555692463648D+00   0.550060331925639D-02   0.526555692463648D+00   0.550060331856020D-02   0.526555692463648D+00   0.550060331419108D-02
 0.430D+01  0.100D+02 -0.637291215759541D+57   0.243857823559396D+57 -0.637291215759540D+57   0.243857823559391D+57 -0.637291215759541D+57   0.243857823559401D+57 -0.637291215759541D+57   0.243857823559391D+57
 0.630D-01  0.120D+02  0.141361484205804D+00   0.499527520944462D+00  0.141361484205804D+00   0.499527520944458D+00  0.141361484205804D+00   0.499527520944463D+00  0.141361484205807D+00   0.499527520944462D+00
 0.630D-01  0.150D+02  0.706036692953861D-03   0.707112019154043D+00  0.706036692961832D-03   0.707112019154043D+00  0.706036692955051D-03   0.707112019154043D+00  0.706036692950090D-03   0.707112019154049D+00
 0.630D-01  0.200D+02  0.416621866088240D+00   0.498381252020572D+00  0.416621866088239D+00   0.498381252020566D+00  0.416621866088240D+00   0.498381252020578D+00  0.416621866088265D+00   0.498381252020613D+00
 0.130D+02  0.160D+02  0.300533459242459D+282 -0.370003657663020D+282  0.300533459242463D+282 -0.370003657663035D+282  0.300533459242482D+282 -0.370003657663059D+282  0.300533459242508D+282 -0.370003657663069D+282
 0.100D+01  0.100D-01  0.780050492928563D+00   0.523753815166502D-06  0.780050492928563D+00   0.523753815262484D-06  0.780050492928563D+00   0.523753815167000D-06  0.780050492928563D+00   0.523753811089223D-06
 0.550D+01  0.100D-03  0.478421381863831D+00  -0.923879994230660D-04  0.478421381863832D+00  -0.923879994230704D-04  0.478421381863831D+00  -0.923879994230661D-04  0.478421381863846D+00  -0.923879994230660D-04
 0.190D+01  0.100D+01  0.362261208092774D+23   0.699962832595258D+24  0.362261208092718D+23   0.699962832595254D+24  0.362261208092757D+23   0.699962832595250D+24  0.362261208093516D+23   0.699962832595250D+24
 0.100D+01  0.280D+02  0.319560919884034D+35  -0.905034298961766D+36  0.319560919883390D+35  -0.905034298961765D+36  0.319560919884066D+35  -0.905034298961759D+36  0.319560919885436D+35  -0.905034298961693D+36
 0.150D+02  0.150D+02  0.512490992884674D+305  0.512490992884674D+305  0.512490992884655D+305  0.512490992884655D+305  0.512490992884682D+305  0.512490992884682D+305  0.512490992884647D+305  0.512490992884655D+305
 0.260D+02  0.000D+00  0.499994235272720D+00   0.000000000000000D+00   0.499994235272720D+00   0.000000000000000D+00   0.499994235272720D+00   0.000000000000000D+00   0.499994235272941D+00   0.000000000000000D+00
 0.000D+00  0.266D+02  0.000000000000000D+00   0.492368098869665D+00   0.000000000000000D+00   0.492368098869665D+00   0.000000000000000D+00   0.492368098869666D+00   0.000000000000000D+00   0.492368098869892D+00
```